\documentclass{elsart}

\usepackage{amsmath}
\usepackage{amssymb}
\usepackage{graphics}
\usepackage{graphicx}
\usepackage{natbib}

\DeclareMathOperator{\Prob}  {\mathsf{P}}
\DeclareMathOperator{\Expec} {\mathsf{E}}

\begin{document}

\begin{frontmatter}


\title{Linear Latent Structure Analysis: 
from Foundations to Algorithms and Applications
\thanksref{thankstogrants}}
\thanks[thankstogrants]{}


\author{I. Akushevich, M. Kovtun, A.I. Yashin and K.G. Manton}

\address{Center for Demographic Studies, Duke University
Durham, North Carolina 27708-0408, USA}

\begin{abstract}
A new statistical technique for constructing 
linear latent structure (LLS) 
models from available data, supported by well established 
theoretical results and an 
efficient algorithm, is presented. The method 
reduces the problem of estimating LLS model parameters to a sequence of 
linear algebra problems.  This assures a low computational complexity and 
an ability to handle large scale data that involve thousands of 
variables. An overall computational scheme and all its components 
are discussed in detail. Simulation experiments demonstrate the excellent 
performance 
of the  algorithm in 
reconstructing model parameters. Step-by-step analysis of a demographic 
survey is presented as an example. 

The technique is useful for the analysis of high-dimensional 
categorical data (e.g., 
demographic surveys, gene expression data) where the detection, evaluation, 
and interpretation of a underlying 
latent structure are required. 

\end{abstract}


\end{frontmatter}

\section{Introduction}
Categorical data often occur in behavior sciences, 
epidemiology, genetics, education, and marketing. A common feature of 
modern categorical data collection (e.g., demographic surveys or gene 
expression data) is their high dimensionality with respect to both 
to sample size and
the numbers of categorical variables measured for each 
individual. 

Limited numbers of categorical variables may not precisely characterize an 
object. This can be due to the stochastic nature of observed 
objects, their heterogeneity, or measurement error. 
Often the phenomenon of interest cannot be observed directly, so a series of 
measures has to be made, each  related to 
latent variables of interest. As a consequence, 
the set of categorically scored measurements is 
strongly related. Researchers analyzing such data have to: 
i) extract the signal from data with a large stochastic 
component, ii) extract a latent (unobservable) structure, 
related to the categorical data, and iii) parameterize the latent 
structure to interpret it substantively. The
identification and 
description of a latent structure are not sufficient when estimation of 
individual parameters are required. Such problems are informally solved by 
a physician making a diagnosis: he/she estimates an unobservable variable of 
interest (presence of a disease) based on categorical 
measures (description of symptoms and answers to physician's questions) 
and estimated population characteristics (physician's professional 
experience). Likewise, having a set of categorical measurements (e.g., a 
demographic survey) made on individuals, a researcher 
would like to discover i) properties of the population and ii) properties of 
the individuals. Methods which attempt to achieve these simultaneously are 
known as  latent structure analysis. 

One of the best known such methods is the latent class model (LCM), which 
can be characterized as a statistical method for finding subtypes of related 
cases (latent classes) from multivariate categorical data (Lazarsfeld and 
Henry, 1968; Goodman, 1974; Clogg, 1995). In LCM, individuals are 
assigned to one of several homogeneous classes. This requires 
estimation of the individual 
latent variable (class number). Other specific models of this family such as 
item response theory and Rasch models, discrete latent class models (Heinen, 
1996), latent distribution analysis (Mislevy, 1984; Uebersax {\&} 
Grove, 1993; Qu, Tan {\&} Kutner, 1996), differ by the assumptions made about
the 
latent variable(s). One method for identifying the latent 
structure in large categorical data sets with simultaneous evaluation of 
individual scores 
in a state space is Grade of Membership (GoM) analysis. GoM was introduced 
in Woodbury and Clive (1974); see Manton et al. (1994) for a detailed 
exposition and additional references. Statistical properties of GoM models 
were analyzed by Tolley and Manton (1992). 

All these models belong to the family 
of  latent structure analysis that considers a number of categorical 
variables measured on each individual in a sample with the intent of 
discovering the
properties of both a population and individuals composing the 
population. Different approaches to this general problem, and recent
developments, are 
described in reports collected in Hagenaars and McCutcheon 
(2002). The relation of specific latent structure models is 
discussed by Uebersax (1997) and Erosheva (2005). Comparative analyses of 
methods of LSA and traditional statistical approaches (Agresti 2002) 
attempting to find the
correspondence of population parameters of interest and sample 
statistics are presented by Powers and Xie (2000).

Methods of latent structure analysis can be 
reformulated in terms of statistical mixtures and mixing distributions. 
Bartholomew (2002) noted the theory of mixing distributions 
as a common mathematical framework of 
latent structure models with abilities ``to clarify, simplify, and unify the
disparate developments spanning over a century.'' 
The main assumption of the latent structure 
analysis is the local independence assumption. 
Being formulated in terms of the theory of mixing distributions,
it states that the observed joint distribution (of
categorical random variables) can be described as a mixture of
independent distributions. The 
mixing distribution is considered as a distribution of the latent variable(s) 
that contains latent information regarding the 
phenomenon under study. Specific models of latent structure analysis 
vary by assumptions regarding properties of the mixing distribution. 

LCM, GoM, and many other methods of latent structure 
analysis use maximum likelihood 
for parameter estimation. Although a part of parameters describing a latent 
structure are regular (structural) parameters, parameters corresponding to 
individual states are so-called ``nuisance'' parameters, the number of which 
may increase 
with sample size. There are a series of non-trivial mathematical tasks which 
have to be solved to estimate such individual parameters. 
First is the problem of identifiability, which is 
especially difficult for large samples and a large number of estimated 
parameters. Second is the problem of consistency, which is very difficult to
prove when  
structural and nuisance parameters are estimated simultaneously. 
Third is the availability and quality of algorithms to 
perform estimations. Often additional assumptions not required by a model 
are made  to facilitate  likelihood maximum. No general 
theorems are available to address these questions, so all of these points have 
to be considered separately for each task involved in dealing with nuisance 
parameters (e.g., Haberman, 1995; Erosheva, 2002). 

Recently Linear Latent Structure (LLS) analysis has 
been proposed to model high dimensional 
categorical data (Kovtun et al., 2005a,b,c).  
LLS has been formulated using mixing distribution theory. Similar to other 
latent structure analyses, the goal of LLS analysis is to derive both the 
properties of a population and individuals using 
discrete measurements. LLS, however, does not use a maximum likelihood 
approach for parameter estimation. Instead, it uses a new estimator, 
where LLS parameter
estimates are 
solutions of a quasilinear system of equations. For this estimator, it is 
possible to prove consistency, to formulate conditions for identifiability, 
and to formulate a high-performance algorithm allowing one to handle 
datasets involving thousands of categorical variables.

The LLS approach is briefly reviewed 
in Section 2. Specifically, the structure of the data base, 
the LLS problem formulations, their basic properties, and two examples are 
discussed. A detailed description of the LLS algorithm and its implementation 
is in Section 3. Attention is paid 
not only to the components of LLS analysis, but also to auxiliary 
algorithms, and to procedures for choosing a basis which may require prior 
information about the study phenomenon. Section 4 is devoted to applications,
and has 
 two parts. The first 
 contains a description of simulation experiments 
designed to check the quality of reconstruction of  model components. 
The second  discusses how the model  is applied to a 
specific dataset. Section 5 summarizes results and includes 
discussion and conclusions. 

\section{Linear Latent Structure Analysis}

\subsection{Structure of datasets and population characteristics. }
The typical dataset analyzed by methods of latent structure analysis can be 
represented by the $I\times J$ matrix constituted by categorical outcomes 
$X_j^i $ of $J$ measurements on $I$ individuals, where $i=1,\ldots ,I$ and 
$j=1,\ldots ,J$ run over individuals and measurements, respectively. Each 
row in the matrix corresponds to an individual and contains an 
individual response pattern, i.e., a sequence of $J$ numbers with the 
$j$th number
running from 1 to the number of responses
$L_j $ for that variable. In most cases $L_j$ ranges from 2 to 5-10, and
rarely exceeds several 
dozens. Thus, the results of a  survey are represented by $I$ 
measurements of random variables $X_1 ,\ldots ,X_J $, with the set of 
outcomes of the $j$th measurement being $\left\{ {1,...,L_j } \right\}$. 
The joint distribution of 
random variables $X_1 ,\ldots ,X_J $ can be described by the elementary 
probabilities,
\begin{equation}
\label{eq1}
p_\ell =\Prob \left( {X_1 =\ell _1 \;\;{\rm and}\;\;\cdots 
\;\;{\rm and}\;\; X_J =\ell _J } 
\right),
\end{equation}
where $\ell =(\ell _1 ,...,\ell _J )$ is an individual response pattern and 
$\ell _j \in \left\{ {1,...,L_j } \right\}$. To include into consideration
marginal probabilities, we allow some components of $\ell$ to be $0$'s. For
example, for three binary variables, 
$$
p_{(2,0,1)}=\Prob (X_1=2\;\;{\rm and}\;\; X_3=1)=p_{(2,1,1)}+p_{(2,2,1)}
$$ 
Values of these probabilities 
$p_\ell $ (and only these) are directly estimable from the 
observations. If $I_\ell $ is the number of individuals with pattern $\ell $, 
consistent estimates for $p_\ell $ are given by frequencies $f_\ell =I_\ell 
/I$. 

\subsection{LLS task: statistical, geometrical, and mixing 
distribution points of view.}

The problem in LLS analysis is to evaluate dimension 
of a hidden space, identify
its location in the space of larger dimension, and to evaluate hidden
individuals' characteristics (coordinates in the latent sub-space) from the
data. The LLS analysis is based on two assumptions. The first is the assumption
about ``local independence'', which is common for all methods of latent
structure analysis. The second is specific for LLS analysis. It is about 
existence of
low-dimensional linear subspace associated with the latent structure.
We present LLS in terms of the theory of mixing
distributions, and then
discuss its specific assumption 
from statistical and geometrical points of view. 

Population characteristics are completely described by the 
joint distribution of 
random variables $X_1 ,\ldots ,X_J $ presented by probabilities (\ref{eq1}). 
Among all possible joint distributions, one can 
distinguish independent distributions, i.e. distributions satisfying,
\begin{equation}
\label{eq3}
p_\ell =\Prob \left( {X_1 =\ell _1 \,\,{\rm and}\,\,\cdots \,{\rm and}\,X_J =\ell _J } 
\right)=\prod\nolimits_j {\Prob \left( {X_j =\ell _j } \right)} .
\end{equation}
The description of an independent distribution law requires only knowing 
$\Prob \left( {X_j =\ell _j } \right)$ denoted below as $\beta _{jl}$. 
Vectors of probabilities 
$\beta=(\beta _{11},\ldots,\beta _{JL_J}) $ belong to vector space $R^{\left| 
L \right|}$, where $\left| L \right|=\sum\nolimits_j {L_j } $. Indexes of the
vector components run over all possible pairs of $jl$, i.e., corresponding to 
probabilities of the first outcome to the first question, of the second 
outcome to the first question, and so on. Requirements for $\beta _{jl} $ 
to be probabilities restricts their domain in the vector space by 
\begin{equation}\label{eq2a}
\sum\limits_{l=1}^{L_j } {\beta _{jl} } =1 
\qquad {\rm and} \qquad \beta _{jl} 
\ge 0. 
\end{equation}
This domain represents the direct product of $J$ unit simplex 
of dimensions $L_j $.

Since variables $X_1 ,\ldots ,X_J $ in general case 
are not independent, the observed
distribution $\{p_\ell \}$ cannot be described by the product 
of independent distributions,
but it can be exactly described  
as a mixture of independent 
distributions. This means that each set of independent probabilities contributes
to observed distribution with a weight function. This weight function is called
mixing distribution. It is defined in the space of independent distributions,
i.e. for each vector of probabilities $\beta$ satisfying (\ref{eq2a}). 
Let $F(\beta) $ be the cumulative distribution function of the mixing
distribution. Probabilities $p_\ell$ are represented as,
\begin{equation}
\label{eq4}
p_\ell =\int dF(\beta) \prod\nolimits_{j=1}^J \beta _{j\ell_j }   
\end{equation}
Thus, latent 
structure analysis searches for a representation of the observed 
distribution as a mixture of independent distributions.

Any distribution 
$\{p_\ell \}$ can be represented as a mixture, 
so representation (\ref{eq4}) 
does not restrict the family of distributions and further specifications 
are required. They are formulated as restrictions on the support of mixing 
distribution or, equivalently, on a set of mixed independent distributions.  
The LLS specific assumption is that this set is restricted to be a 
$K$-dimensional linear subspace of the space 
of independent distributions, i.e.,  the mixing 
distribution is supported by the 
linear subspace spanned by $K$ basis vectors $\lambda ^1,\ldots ,\lambda ^K$.
Below this LLS assumption is considered from the point of view of pure 
statistical analysis and the geometry of the task.

Individual characteristics are described by individual
probabilities $\beta _{jl}^i =\Prob (X_j^i 
=l)$ of specific outcomes ($i=1,\ldots ,I$ runs over individuals).

  The LLS assumption about the existence of a low-dimensional linear 
subspace supporting the mixing distribution is essentially equivalent to the 
assumption that there exists a $K$-dimensional random vector $G$ such that 
for every $j$ a regression of $Y_{jl} $ (random variable $Y_{jl} $ 
equaling 1 if $X_j =l$ and 0 otherwise) on $G$ is linear.
The regression 
equation relates the expectation of $Y_{jl} $, which is $\beta _{jl} $, 
to  the random vector $G$. 
If a specific value of $G$ is associated
 with individual $i$ (so-called  
LLS scores $g_{ik} $), then the regression takes the form, 
\begin{equation}
\label{eq2}
\beta _{jl}^i =\sum\limits_{k=1}^K {g_{ik} \lambda _{jl}^k } .
\end{equation}
The sense of the regression coefficients $\lambda _{jl}^k $ and model 
restrictions is clarified by analyzing the geometry of the LLS 
task. 

Vectors of individual probabilities $\beta^i=\{ \beta^i_{jl}\}$, of
individual responses $Y^i=\{Y_{jl}^i\} $ and the regression coefficients 
$\lambda^k=\{\lambda _{jl}^k \}$ lie in the permitted domain (\ref{eq2a}) of the
space of independent distributions. From a 
geometric point of view, LLS searches a $K$-dimensional subspace (represented 
by $\lambda _{jl}^k )$ in this space, which is the 
``closest'' to the set of $I$ points representing individual outcomes $Y_{jl}^i $. 
This linear subspace is defined by its basis $\lambda ^1,\ldots ,\lambda ^K$, 
 so to find the $K$-dimensional subspace means 
finding a basis, $\lambda _{jl}^k $, ($k=1,\ldots ,K)$. Every basis 
$\lambda ^1,\ldots ,\lambda ^K$ defines a family of regression coefficients and 
vice versa. The 
complete set of restrictions in the LLS task allowing to consider $\beta 
_{jl}^i $ and $\lambda _{jl}^k $ as probabilities, is,
\begin{equation}\label{eq2b}
\sum\limits_{l=1}^{L_J } {\lambda _{jl}^k } =1,
\quad
\lambda _{jl}^k \ge 0,
\quad
\sum\limits_{k=1}^K {g_{ik}} =1 
\quad {\rm and} \quad 
\sum\limits_k {g_{ik} } 
\lambda_{jl}^k \ge 0.
\end{equation}
LLS scores $g_{ik}$
characterizing an individual $i$ are then estimated as the 
expectation of vector $G$, conditional on the respondent's answers. Basis 
vectors of the subspace can be interpreted as probabilities and can define 
the ``pure types'' (Manton et al., 1994). In this sense, the model decomposition 
(\ref{eq2}) has the 
interpretation of a decomposition over pure types or over ``ideal persons'' 
whose individual probabilities are basis vectors of the subspace.

Summarizing, one can say that the LLS model approximates the observed distribution
of $X_1 ,\ldots ,X_J $ by a mixture of independent distributions with a mixing
distribution supported by a $K$-dimensional subspace of the space of independent
distribution. To specify such a model distribution it is sufficient to define the
following LLS parameters:
\begin{enumerate}
\item A basis $\lambda ^1,\ldots ,\lambda ^K$ of the space that supports the
mixing distribution.
\item Conditional moments $\Expec (G|X=\ell)$.
\end{enumerate}  
This set of model parameters is not the only set possible. We chose 
it because of a
number of useful properties listed below.

{\bf Property 1}. The mixing distribution can be estimated in the style of an 
empirical distribution, i.e., when the estimator is a distribution concentrated
in points $\Expec (G|X=\ell)$ with weights $f_\ell$.

{\bf Property 2}. The conditional expectations $\Expec (G|X=\ell)$ provide
knowledge about individuals. These conditional expectations can be considered as
coordinates in a phase space, to which all individuals belong. The ability to
discover the phase space and determine individual positions in it is a valuable
feature of LLS analysis.   

{\bf Property 3}. When the number of measurement, $J$, tends to infinity, the
individual conditional expectations $g_i=\Expec (G|X=\ell^{(i)})$, 
where $\ell^{(i)}$ is the vector of responses of individual $i$, converge to the
true value of the latent variable for this individual, and estimates of the
mixing distribution converge to the true one, provided some regularity
conditions (Kovtun et al., 2005d).

\subsection{Moment matrix and the main system of equations. }

Parameter estimation is based on two  facts (Kovtun et al., 2005a,b) 
 formulated in terms of the conditional and unconditional moments of the
mixing distribution. The 
first  is that columns of moment matrix belong exactly to the desired linear 
space. The second is that they obey the main system of equations.

\subsubsection{Unconditional moments and the moment matrix}

The first set of values in which we are interested consists of 
 the unconditional moments of the mixing distribution 
$F(\beta )$,
\begin{equation}
\label{eq5}
 M_\ell =\int dF(\beta) \prod\nolimits_{j:\ell _j \ne 0} \beta_{j\ell _j }  
p_\ell 
\end{equation}
Note an important fact regarding the above equation. The value on the
left-hand-side, $M_\ell$, is a moment of {\it mixing distribution}, while the
value on the right-hand-side, $p_\ell$, comes from the {\it joint distribution of 
$X_1 ,\ldots ,X_J $}; the equality of these values is a direct corollary of the
definition of mixture. The existence of their connection between two distinct
distributions is crucial for LLS analysis. 

The first corollary of eq. (\ref{eq5}) is that the unconditional 
moments are directly estimable from  data and, therefore, the 
frequencies $f_\ell $ of response patterns $\ell $ observed in a sample are 
consistent and efficient estimators for conditional moments $M_\ell$. 

Recall that we allow some components of response pattern $\ell$ to be 0.
In this case $p_\ell$ are marginal probabilities. In the definition of
$M_\ell$ the multipliers, corresponding to 0 components of $\ell$, are excluded from
the product. Thus, the order of moment $M_\ell$ is equal to the number of non-zero
components in $\ell$.   

All moments defined in (\ref{eq5}) are estimable by frequencies; however, this 
definition does not cover all moments of a certain order. For example, moments
of second order with 
$\beta _{jl_1} $ and $\beta _{jl_2} $, (i.e., with the same $j$) 
are not estimable. 
This  arises because the data do not include double answers to the same 
question. One can notice that i) all moments of first order are estimable, ii)
the proportion of estimable moments decreases with the increase of order, and iii)
no moments of order $J+1$ and higher are estimable.  

The moment matrix is constructed from moments of order up to $J$ using the
following formal rules:
\begin{enumerate}
\item Rows of the moment matrix are indexed by response patterns containing
exactly one non-zero component or, equivalently, by pair indexes $jl$. Thus, the
moment matrix contains $|L|$ rows, and their columns can be considered as
vectors in $R^{|L|}$.    
\item Columns of the moment matrix are indexed by all possible response patterns,
including a response pattern containing all 0's. The first column is indexed by
response pattern $(0,\ldots,0)$; the next $|L|$ columns are indexed by response
patterns containing one non-zero component, and so on. 
\item The element on the intersection of row $\ell '$ and column $\ell ''$ is
$M_{\ell'+\ell''}$, if $\ell ''$ has 0 at the position of the only non-zero
component of $\ell '$ (in this case, $\ell '+\ell''$ is a meaningful response
pattern; otherwise, the question mark is placed on the position of intersection
of row $\ell '$ and column $\ell ''$). For example, the element of the moment
matrix in row $(1,0,0)$ and column $(0,2,2)$ is $M_{1,2,2}$, and element in row
$(1,0,0)$ and column $(1,0,2)$ is a question mark.    
\end{enumerate}

Equation (\ref{eq7}) gives an example of a portion of the moment matrix for the case 
of $J=3$ dichotomous variables, i.e., $L_1 =L_2 =L_3 =2$.
\begin{equation}
\label{eq7}
\left(
\begin{array}{ccccccccc}
M_{(100)}   &
    ?           &  ?            &
    M_{(110)}   &  M_{(120)}    &
    M_{(101)}   &  M_{(102)}    &
    ?           &  \cdots \phantom{\vdots}       \\
M_{(200)}   &
    ?           &  ?            &
    M_{(210)}   &  M_{(220)}    &
    M_{(201)}   &  M_{(202)}    &
    ?           &  \cdots \phantom{\vdots}       \\
M_{(010)}   &
    M_{(110)}   &  M_{(210)}    &
    ?           &  ?            &
    M_{(011)}   &  M_{(012)}    &
    ?           &  \cdots \phantom{\vdots}       \\
M_{(020)}   &
    M_{(120)}   &  M_{(220)}    &
    ?           &  ?            &
    M_{(021)}   &  M_{(022)}    &
    ?           &  \cdots \phantom{\vdots}       \\
M_{(001)}   &
    M_{(101)}   &  M_{(201)}    &
    M_{(011)}   &  M_{(021)}    &
    ?           &  ?            &
    M_{(111)}   &  \cdots \phantom{\vdots}       \\
M_{(002)}   &
    M_{(102)}   &  M_{(202)}    &
    M_{(012)}   &  M_{(022)}    &
    ?           &  ?            &
    M_{(112)}   &  \cdots \phantom{\vdots}      
\end{array}
\right)
\end{equation}
In this example, places for 
inestimable moments are filled by question marks. The first column of the 
moment matrix contains moments of the first order, when only one specific 
outcome of one specific question is taken into account. There are no 
inestimable moments in the first column. Elements of this column
 can be denoted
as components of vectors in $R^{|L|}$, i.e., as $M_{jl}$. 
The next six ($\vert L\vert $ in general) columns correspond to second-order moments. 
Blocks of diagonal 
elements are not estimable. Second-order moments can be also denoted via pair $jl$
of indexes as $M_{jl;j'l'}$. The last shown column corresponds to third order 
moments. The notation $M_{jl}$ and $M_{jl;j'l'}$ is used below for specific
columns of the moment matrix. 

The part of the moment matrix consisting of second-order moments (which is
$|L|\times |L|$ square matrix) together with the column of first-order moments
is of special interest. 
A well-know fact is that if a distribution in an $n$-dimensional Euclidean space is
carried by a $k$-dimensional linear manifold, then the rank of the 
covariance matrix
is equal to $k$, and the position of the manifold can be derived from the 
covariance
matrix. This fact is the cornerstone of principal component analysis. Our method is
based on similar ideas, adapted to having an incomplete set of
second-order moments.   
For a small $J$ (as in the example), there is a relatively large 
fraction of non-estimable components in the second-order part of the 
moment matrix. For increasing $J$, this 
fraction rapidly decreases. 

For a moment matrix $M$ let its completion $\bar {M}$ be a matrix obtained 
from $M$ by replacing question marks with arbitrary numbers. The main fact 
with respect to the moment matrix is that the moment matrix always has a 
completion in which all columns belong to the supporting plane $\Lambda $. 
Thus, if the estimable part of the
moment matrix has sufficient rank (which is the case in 
non­degenerate situations,) a basis in $\Lambda $ may be obtained from this 
matrix. As we have a consistent estimator of the moment matrix in the form of a 
frequency matrix, the supporting plane may be consistently estimated. 

The mixing distribution in the LLS model is supported by the intersection of the 
linear 
subspace with the polyhedron defined by (\ref{eq2a}). 
This intersection $P_\beta $ can be 
parameterized by $g$'s in accordance with (\ref{eq2}). In Appendix, we
demonstrate that integrals in (\ref{eq5}) over $P_\beta$ can be reduced to
integrals in the linear subspace over the polyhedron $P_g$, and how distribution
functions $F(\beta)$ and $F(g)$ relate to each other. 
Then unconditional moments $M_\ell$ can be formally represented as
\begin{equation}
\label{eq8a}
M_\ell=p_\ell=\int dF(g)
\prod\nolimits_{j:\ell _j \ne 0} \sum\nolimits_k g_k \lambda^k_{j\ell _j }
\end{equation}

Representation (\ref{eq8a}) helps to illustrate the main properties of the moment
matrix: i) the rank of the matrix is $K$ and ii) all columns of the matrix
belong to a linear subspace with a basis represented by 
vectors $\lambda^k$. Indeed, since $jl$ indexes running in each column of the
moment matrix are only
involved in $\lambda$'s, which are independent of integration 
variable ${\bf g}$, and,
therefore, which can be removed from the integral. The remaining integration 
over ${\bf g}$ is common for all elements of the moment matrix of a 
certain order.
The first column and second order
matrix are represented,
\begin{eqnarray}\label{eq9aa}
M_{jl}&=&\int dF(g)
\sum\limits_k g_k \lambda^k_{jl} 
=\sum\limits_k \left( \int dF(g) g_k \right) \cdot \lambda^k_{jl}
=\sum_k I_k \lambda^k_{jl}
\nonumber
\\
M_{jl;j'l'}&=&
\int dF(g)
\left[\sum\limits_k g_k \lambda^k_{jl}\right] 
\left[\sum\limits_{k'} g_{k'} \lambda^{k'}_{j'l'}\right] 
=\sum\limits_{kk'}\lambda^k_{jl}\lambda^{k'}_{j'l'}
\int dF(g)
g_k  g_{k'} \nonumber\\&=&
\sum\limits_k \left( \sum\limits_{k'} \lambda^{k'}_{j'l'}
\int dF(g) g_k g_{k'} \right) \cdot \lambda^k_{jl}
= \sum_{k} I_{k j'l'} \lambda^k_{jl}
\end{eqnarray}
where $I_k$ and $I_{k j'l'}$ are some numeric coefficients. From 
(\ref{eq9aa}) we see that each column of moment matrix can be represented as a 
linear combination of $K$  vectors $\lambda^k$.

\subsubsection{Unconditional moments and main system of equations}

Another set of the values of interest are the conditional moments $\Expec(G_k \vert X=\ell )$, 
which express  knowledge of the state of individuals based on measurements. 
They are not directly estimable from  
observations. The goal of LLS analysis is to obtain  estimates for these 
conditional moments. Explicit expressions for those of the lowest order are 
obtained using the Bayes theorem (Kovtun et al. 2005a),
\begin{equation}
\label{eq9}
\Expec(G_n \vert X=\ell )=\int dF(g){g_n } \frac{\prod\nolimits_{j:\ell _j \ne 0} 
 \; {\sum\nolimits_k {g_k \lambda _{jl}^k } } }{M_\ell (\mu _\beta )}
\end{equation}
Analogously, higher conditional moments, including products of components of 
$G$, can be constructed. As can be seen explicitly from (\ref{eq8}) and (\ref{eq9}), 
the relation 
of conditional and unconditional moments in LLS analysis can be described 
as,
\begin{equation}
\label{eq10}
\sum\nolimits_k {\lambda _{jl}^k } \cdot \Expec(G_k \vert X=\ell )\frac{M_{\ell +l_j }}{M_\ell},
\end{equation}
where vector $\ell $ contains 0 in position $j$, and $\ell +l_j $ contains $l$ 
in this position. Similar equations can be written for conditional moments 
of higher orders. We refer to the system of equations relating conditional and 
unconditional moments as the main system of the equation. Kovtun et al, (2005a) 
proved the following properties of solutions of the main system of 
equations: i) any basis $\lambda _{jl}^k $ of $\Lambda $ together with 
conditional moments $\Expec(G_k \vert X=\ell )$ calculated on this basis give a 
solution of the main system of equation; and ii) under regular conditions, 
every solution of the main system of equations gives a basis of $\Lambda $ 
and conditional moments calculated in this basis. Note, that equation 
(\ref{eq10}) is linear with respect to conditional moments.

The described properties of the moment matrix and solutions of the main 
system of equations suggest an efficient algorithm to obtain LLS estimates. 
First, a basis of the supporting plane can be obtained from the moment 
matrix, and second, 
conditional moments can be found by solving a linear system of equations. 
An important property of conditional moments is that (Kovtun et al., 
2005d): if a population is represented by a set of true LLS scores, then 
unconditional moments give good approximation for them for finite $J$, and 
converge to true LLS scores when $J\to \infty $. 

\subsection{Two illustrative examples.}

Before going into detail for the algorithm and to realistic tasks of data 
analysis, we consider two simple illustrative examples. For both of them, 
assume $K=2$, three dichotomous variables ($J=3)$, and the basis vectors are 
$\lambda ^1=(1,0;1,0;1,0)$ and $\lambda ^2=(\raise0.5ex\hbox{$\scriptstyle 
1$}\kern-0.1em/\kern-0.15em\lower0.25ex\hbox{$\scriptstyle 
2$},\raise0.5ex\hbox{$\scriptstyle 
1$}\kern-0.1em/\kern-0.15em\lower0.25ex\hbox{$\scriptstyle 2$};0,1;0,1)$. 
Then the independent distributions being mixed are defined by vectors:
\begin{equation}
\beta=g_1\lambda^1+g_2\lambda^2=g_1\lambda^1+(1-g_1)\lambda^2, \qquad 
0\leq g_1\leq 1.
\end{equation} 
Thus, a mixing distribution can be given one dimensional p.d.f. $\rho(g_1)$. 
 For the first task, we assume that 
the mixing distribution is uniform ($\rho (g_1)=1\cdot \theta(g_1)\cdot
\theta(1-g_1)$. In the second case we 
assume the mixing distribution is concentrated at two points with $g_1 =0.1$ 
and $g_1 =0.4$ ($\rho (g_1 )=\raise0.5ex\hbox{$\scriptstyle 
1$}\kern-0.1em/\kern-0.15em\lower0.25ex\hbox{$\scriptstyle 2$}[\delta (g_1 
-\raise0.5ex\hbox{$\scriptstyle 
1$}\kern-0.1em/\kern-0.15em\lower0.25ex\hbox{$\scriptstyle {10}$})+\delta 
(g_1 -\raise0.5ex\hbox{$\scriptstyle 
2$}\kern-0.1em/\kern-0.15em\lower0.25ex\hbox{$\scriptstyle 5$})])$. 
Unconditional moments are calculated using (\ref{eq5}). Moment matrices for both 
cases are
\begin{equation}
\label{eq11}
\begin{array}{ccc}
\left( 
\begin{array}{*{7}c}
{\textstyle{3 \over 4}}  & 
{\textstyle{\it 7 \over {12}}}  & 
{\textstyle{\it 1 \over 6}}  & 
{\textstyle{1 \over 2}}  & 
{\textstyle{1 \over 4}}  & 
{\textstyle{5 \over {12}}}  & 
{\textstyle{1 \over 3}}  
\\
 {\textstyle{1 \over 4}}  & 
 {\textstyle{\it 1 \over 6}}  & 
{\textstyle{\it 1 \over 12}}  & 
{\textstyle{1 \over 8}}  & 
{\textstyle{1 \over 8}}  & 
{\textstyle{1 \over {12}}}  & 
{\textstyle{1 \over 6}}  
\\
 {\textstyle{5 \over 8}}  & 
 {\textstyle{1 \over 2}}  & 
{\textstyle{1 \over 8}}  & 
{\textstyle{\it 7 \over {16}}}  & 
{\textstyle{\it 3 \over {16}}}  & 
{\textstyle{3 \over 8}}  & 
{\textstyle{1 \over 4}}  
\\
 {\textstyle{3 \over 8}}  & 
 {\textstyle{1 \over 4}}  & 
{\textstyle{1 \over 8}}  & 
{\textstyle{\it 3 \over {16}}}  & 
{\textstyle{\it 3 \over {16}}}  & 
{\textstyle{1 \over 8}}  & 
{\textstyle{1 \over 4}}  
\\
 {\textstyle{{\rm {\bf 1}} \over {\rm {\bf 2}}}}  & 
 {\textstyle{{\rm {\bf 5}} \over {{\rm {\bf 12}}}}}  & 
 {\textstyle{1 \over {12}}}  & 
 {\textstyle{3 \over 8}}  & 
 {\textstyle{1 \over 8}}  & 
{\textstyle{\it 1 \over 3}}  & 
{\textstyle{\it 1 \over 6}}  
\\
 {\textstyle{1 \over 2}}  & 
 {\textstyle{1 \over 3}}  & 
{\textstyle{1 \over 6}}  & 
{\textstyle{1 \over 4}}  & 
{\textstyle{1 \over 4}}  & 
{\textstyle{\it 1 \over 6}}  & 
{\textstyle{\it 1 \over 3}}  
\\
\end{array} 
\right)
& \,\,\,\,\mbox{and}\,\,\,\ &
\left( 
\begin{array}{*{7}c}
{\textstyle{5 \over 8}}  & 
{\textstyle{\it {317} \over {800}}}  & 
{\textstyle{\it {183} \over {800}}}  & 
{\textstyle{{451} \over {1600}}}  & 
{\textstyle{{549} \over {1600}}}  & 
{\textstyle{{67} \over {400}}}  & 
{\textstyle{{183} \over {800}}}  
\\
{\textstyle{3 \over 8}}  & 
{\textstyle{\it {183} \over {800}}}  & 
{\textstyle{\it {117} \over {800}}}  & 
{\textstyle{{249} \over {160}}}  & 
{\textstyle{{351} \over {1600}}}  & 
{\textstyle{{33} \over {400}}}  & 
{\textstyle{{117} \over {400}}}  
\\
 {\textstyle{7 \over {16}}}  & 
 {\textstyle{{451} \over {1600}}}  & 
 {\textstyle{{249} \over {160}}}  & 
 {\textstyle{\it {653} \over {3200}}}  & 
 {\textstyle{\it {747} \over {3200}}}  & 
 {\textstyle{{101} \over {800}}}  & 
 {\textstyle{{249} \over {800}}}  
 \\
 {\textstyle{9 \over {16}}}  & 
 {\textstyle{549 \over {1600}}}  & 
 {\textstyle{{351} \over 1600}}  & 
 {\textstyle{\it {{747}} \over {{3200}}}}  & 
 {\textstyle{\it {{1053}} \over {{3200}}}}  & 
{\textstyle{{{99}} \over {{800}}}}  & 
{\textstyle{{{351}} \over {{800}}}}  
\\
 {\textstyle{{\rm {\bf 1}} \over {\rm {\bf 4}}}}  & 
 {\textstyle{{{\rm {\bf 67}}} \over {{\rm {\bf 400}}}}}  & 
 {\textstyle{{33} \over {400}}}  & 
 {\textstyle{{101} \over {800}}}  & 
 {\textstyle{{99} \over {800}}}  & 
 {\textstyle{\it {17} \over {200}}}  & 
{\textstyle{\it {33} \over {200}}}  
\\
 {\textstyle{3 \over 4}}  & 
 {\textstyle{{183} \over {400}}}  & 
{\textstyle{{117} \over {400}}}  & 
{\textstyle{{249} \over {800}}}  & 
{\textstyle{{351} \over {800}}}  & 
{\textstyle{\it {33} \over {200}}}  & 
{\textstyle{\it {117} \over {200}}}  
\\
\end{array}
  \right)
\end{array}
\end{equation}
Since these matrices were constructed from mixing distributions known a priori,
diagonal blocks in the sub-matrix of the second order
are calculable (marked by the italic font in (\ref{eq11})).
As one can see, the rank of both these matrices is 2. Conditional moments are 
calculated for an outcome pattern. Choose $\ell =(001)$ and $\ell +l_1 
=(101)$. Using (\ref{eq9}) we have,
\begin{equation}
\label{eq12}
\Expec(G_1 \vert X=(001))=2 \mathord{\left/ {\vphantom {2 3}} \right. 
\kern-\nulldelimiterspace} 3\,\,\,\mbox{and}\,\,\Expec(G_2 \vert X=(001))=1 
\mathord{\left/ {\vphantom {1 3}} \right. \kern-\nulldelimiterspace} 3
\end{equation}
for the first example and, 
\begin{equation}
\label{eq13}
\Expec(G_1 \vert X=(001))={17} \mathord{\left/ {\vphantom {{17} {50}}} \right. 
\kern-\nulldelimiterspace} {50}\,\,\,\mbox{and}\,\,\Expec(G_2 \vert X=(001))={33} 
\mathord{\left/ {\vphantom {{33} {50}}} \right. \kern-\nulldelimiterspace} 
{50}
\end{equation}
for the second. Using corresponding elements of $M_\ell $ in (\ref{eq11}) (marked by 
bold text) we can see that l.h.s. and r.h.s of eq. (\ref{eq10}) equal to $5/6$ 
 for first example and $67/100$ for the second:
\begin{equation}
1\cdot {2 \over 3} + {1 \over 2}\cdot  {1 \over 3} = {5/12 \over 1/2}
\qquad {\rm and} \qquad
1\cdot {17 \over 50} + {1 \over 2}\cdot  {33 \over 50} = {67/400 \over 1/4}
\end{equation}
External indexes in this example are $j=1$ and $l=1$. 
 
\section{Algorithms of Linear Latent Structure Analysis}

Parameter estimations in LLS models are 
based on properties of the moment matrix and the main system of equations. 
These
properties allow us to 
reduce a problem of estimating model parameters to a sequence  
of linear algebra problems.
The algorithm based on linear algebra methods 
assures a low computational complexity. 

Data to be 
analyzed are represented by a set of measurements $X_j^i $ (See section 2.1). 
Finding a linear space and individual LLS scores is required. Estimation of 
the model includes four steps: i) estimating the rank of the frequency matrix, 
ii) finding the supporting plane, iii) choosing a basis in the found plane, 
iv) calculating individual conditional expectations and estimating  
mixing distribution. The second and fourth steps are the essence of LLS
parameter estimation problem. The first step is defined as separate 
because 
sometimes the desired dimensionality of the LLS model may be provided by
a researcher, and this step may be skipped. The third step requires
using prior 
information about the processes studied, so it is also examined
separately. 

\subsection{Moment matrix calculation}
An important preceding step that deserves special 
attention is the moment matrix calculation.
The elements of the moment matrix given by $M_\ell $ are approximated by 
observable frequencies defined as $f_\ell =I_\ell /I$, where $I_\ell $ is the 
number of individuals with outcome pattern $\ell $, and $I$ is the total 
number of individuals having certain (not missing) 
outcomes for nonzero elements in $\ell$. 
Columns of a different order have different 
normalizations, e.g., the sum of first-order moments corresponding to question 
$j$ is one (e.g., $M_{(010)} +M_{(020)} =1)$, while sums of columns for this 
$j$ of the second-order sub-matrix are equal to corresponding 
first-order moment (e.g., $M_{(110)} 
+M_{(120)} =M_{(100)} )$. General conditions of summations of the second order
moments written in terms of notation defined after eq. (\ref{eq7}) are,
\begin{equation}  
\sum\limits_{l'=1}^{L_{j'}}M_{jl;j'l'}=M_{jl}.
\end{equation}  

Because of missing data, the property of normalization can be violated. This property,
with or without the renormalization making the sums equaling to one, is required
for the analysis. The
renormalization could provide the property in the case of presence of missing
data, however, this approximation can be
true only assuming missing data are random. More detailed discussion of
procedures to deal with missing data are discussed in Section \ref{sectmissing}.

In addition, a matrix containing standard errors (or 
confidence intervals) of estimates of frequencies is calculated for each 
element of the frequency matrix. Standard errors for binomial distribution, i.e. 
$\sigma _\ell =\sqrt {f_\ell (1-f_\ell )/I_\ell } $, require generalization 
for patterns with small $I_\ell $ as discussed in Brown et al. (2001). A 
generalization based on Wilson's approach (Brown et al., 2001) uses, 
\begin{equation}
\label{eq14}
CI_W =\frac{I_\ell f_\ell +\textstyle{1 \over 2}z_{\alpha /2}^2 }{I_\ell 
+z_{\alpha /2}^2 }\pm \frac{z_{\alpha /2} \sqrt {I_\ell } }{I_\ell 
+z_{\alpha /2}^2 }\sqrt {f_\ell (1-f_\ell )+\frac{z_{\alpha /2}^2 }{4I_\ell 
}} .
\end{equation}
Eq. (\ref{eq14}) recovers the standard Wald's estimates of CIs, 
i.e., $CI_s =p_c \pm 
z_{\alpha /2} \sigma _E $; $z_{\alpha /2} =\Phi ^{-1}(1-\alpha /2)$, where 
$\Phi (x)$ is the standard normal distribution function and $\alpha $ 
denotes the confidence level. 

\subsection{Computational rank of the frequency matrix} 

The frequency matrix can be presented as a sum of the moment matrix with 
rank $K$ and a matrix with a stochastic component. To define the 
dimensionality of the LLS problem, we have to estimate the rank of the 
frequency matrix eliminating the stochastic component. Specifically, we take the 
greatest minor of the frequency matrix that does not contain question marks. 
Then we calculate the singular value decomposition (SVD) and take $K$ equal to 
the number of singular values that are greater than a maximum of the total 
standard deviation estimated as the quadratic sum of standard errors of 
frequencies involved in the minor. 

The  choice of a minor does not essentially influence the 
computational rank of the frequency matrix. Indeed, the geometrically 
specific choice of a minor (e.g. a $n$-dimensional minor of maximal size in 
left low corner of moment matrix) corresponds to projection of a part of 
vectors onto n-dimensional linear subspace. If the real rank of the moment 
matrix is much less than $n$,  it is clear that the rank of the 
projections does not change. 

\subsection{Finding the supporting plane} 

All columns of the moment matrix belong to the 
supporting plane, and as the frequency matrix is an approximation of the 
moment matrix, a natural way to search for the supporting plane is to 
search for a plane that minimizes the sum of distances from it to the 
columns of the frequency matrix. In our case, however, this way is 
complicated by: (a) the frequency matrix is 
incomplete; (b) the statistical inaccuracy of approximation of moments 
$M_\ell $ by frequencies $f_\ell $ varies considerably over elements of 
frequency matrix; and (c) a sought basis should exactly satisfy conditions 
$\sum\nolimits_{l=1}^{L_J } {\lambda _{jl}^k } =1$ for every $k$ and $j$. 
The current prototype of the code overcomes these obstacles by using some 
heuristic methods: (a) An iterative procedure for completion of the frequency 
matrix is used: after a basis of supporting plane is obtained, it is used to 
recalculate completion of the frequency matrix. A new frequency matrix is used 
for adjusting basis calculation etc. (b) Only the first and second order 
moments are examined, so statistical errors of 
different columns in this matrix are compatible. (c) Rotation of each 
simplex (corresponding to each question) to the hyperplane to eliminate one 
degree of freedom. Rotation, but not a simple projection, is required to 
provide the same distances between points in a simplex. 
Items (a) and (c) require explicit consideration.

\subsubsection{Completion of the moment matrix}
We consider the second-order moment 
matrix where for every $\bar {j}$ there are undefined elements corresponding 
to repeated answers to the same question. The intent of completion procedure is to 
approximate these elements, assuming that the supporting subspace $\Lambda $ 
is found. Since only the completed frequency matrix is used for finding  
subspace $\Lambda $, and since the completion procedure uses a basis in the 
sought subspace $\Lambda $, it can be done within the iteration procedure. For 
one iteration step, it is required to find a symmetric matrix $B_{\bar 
{j}} $ of $L_{\bar {j}} \times L_{\bar {j}} $-dimension with positive 
elements and the required summation conditions such that the sum of elements in a column 
(or in a row)  equals to the corresponding moment of the first order, i.e., 
$\sum\nolimits_l {B_{\bar {j},ll'} =M_{\bar {j}l'} } $. Since we know first- and
second-order frequencies ($f_{jl}$ and $f_{jl,j'l'}$; $j\neq j'$), which only
approximate exact moments ($M_{jl}$ and $M_{jl,j'l'}$), special efforts are
required to process the properties of $B_{\bar j}$. Columns of the 
second-order sub-matrix corresponding to question $\bar j$ are 
presented using known frequencies $f_{jl,\bar{j}\bar{l}}$; $j\neq \bar j$
and inestimable elements $B_{\bar {j},l\bar l}$,
\begin{equation}
\left(
\begin{array}{ccc}
\displaystyle f_{11;\bar{j}1} &
\displaystyle \ldots &
\displaystyle f_{11;\bar{j}L_{\bar j}} \\
\displaystyle \ldots &
\displaystyle \ldots &
\displaystyle \ldots \\
\displaystyle f_{1L_1;\bar{j}1} &
\displaystyle \ldots &
\displaystyle f_{1L_1;\bar{j}L_{\bar j}} \\
\displaystyle \ldots &
\displaystyle \ldots &
\displaystyle \ldots \\
\displaystyle B_{{\bar j},11} &
\displaystyle \ldots &
\displaystyle B_{{\bar j},1L_{\bar j}} \\
\displaystyle \ldots &
\displaystyle \ldots &
\displaystyle \ldots \\
\displaystyle B_{{\bar j};L_{\bar j}1} &
\displaystyle \ldots &
\displaystyle B_{{\bar j};L_{\bar j}L_{\bar j}} \\
\displaystyle \ldots &
\displaystyle \ldots &
\displaystyle \ldots \\
\displaystyle f_{11;\bar{j}1} &
\displaystyle \ldots &
\displaystyle f_{J1;\bar{j}L_{\bar j}} \\
\displaystyle \ldots &
\displaystyle \ldots &
\displaystyle \ldots \\
\displaystyle f_{JL_J;\bar{j}1} &
\displaystyle \ldots &
\displaystyle f_{JL_J;\bar{j}L_{\bar j}} 
\end{array}
\right)
\end{equation} 

The completion 
procedure is based on the fact that the rank of the moment matrix is $K$, which 
is much smaller than the dimension of matrix $\vert L\vert $. Therefore, only 
$K$ columns are linearly independent. 
Each column of the moment matrix, being a vector in $K$-dimensional vector 
space, can be expanded over basis vectors $\lambda ^1,\ldots ,\lambda ^K$ 
available after finding the subspace $\Lambda $. Known elements 
$f_{jl;\bar j\bar l}$ (${\bar l} =1,\ldots ,L_{\bar {j}} $ and $j\ne \bar {j})$ 
of 
columns of the moment matrix corresponding to question $\bar {j}$ are 
expanded,
\begin{equation}
f_{jl;\bar j\bar l} =\sum\limits_k {C_k^{\bar j\bar l} } 
\lambda _{jl}^k  \qquad (j\ne \bar 
{j}).
\end{equation}
If coefficients $C_k^{\bar j\bar l} $ are found, matrix $B_{\bar j}$ can be 
constructed as $B_{{\bar j},{\bar l}'{\bar l}} =\sum\limits_k 
C_k^{\bar j\bar l} \lambda ^k_{\bar j{\bar l}'}$, 
The number of known 
components of a vector $f_{jl;\bar j\bar l}$ is greater than the number of basis 
vectors, so coefficients $C_k^{\bar j\bar l}$ can be calculated by ordinary least 
squares methods minimizing residuals. Matrix $B_{\bar {j}} $ has to be 
obtained exactly symmetric, using minimization with restrictions. Lagrangian $Lg_{\bar {j}} $ are written for each  
$\bar j $,
\begin{equation}
\label{eq15}
\begin{array}{l}
 Lg_{\bar {j}} =\sum\limits_{jl:j\ne \bar {j}} {\left( {f_{jl;\bar j\bar l} 
-\sum\limits_k {C_k^{\bar j\bar l}} \lambda _{jl}^k } \right)} 
^2+\sum\limits_{{\bar l} ,k} {\bar\rho _{{\bar l} k} } C_k^{\bar j \bar l} + \\ 
 \,\,\,\,\,\,\,\,\,\,\,\,\,\,\,\,\,\,\,\,\,\,\,\,\,+\sum\limits_{\bar l} 
 {\rho _{\bar l}} 
 \left( {\sum\limits_k {C_k^{\bar j \bar l} -1} } \right)
 -
 \sum\limits_{{\bar l} \ne {\bar l}' } 
 \rho _{{\bar l}{\bar l}'} \sum\limits_k 
 \left( 
 C_k^{\bar j \bar l} \lambda^k_{\bar j {\bar l}'} 
-C_k^{\bar j {\bar l}'} \lambda^k_{\bar j \bar l} 
 \right)
 \\ 
 \end{array}
\end{equation}
where $\rho_{\bar l}  $ and $\rho_{{\bar l}{\bar l}'} $ 
are Lagrange multipliers 
for equality conditions, and $\rho _{{\bar l} k} $ --- for 
inequality conditions. The optimization task to 
find $C_k^\alpha $ with a quadratic minimizing function and 
linear equality ,and inequality conditions are a quadratic programming problem. This 
problem is used by different parts of the 
algorithm (see eqs. (\ref{eq20}) and (\ref{eq23})), so Section 3.7 is devoted to 
this calculation. 

\subsubsection{Removing restrictions} 
The restrictions 
$\sum\nolimits_{l=1}^{L_J } {\lambda _{jl}^k } =1$ are removed by reducing 
the number of rows by $J$ (one for every group of indexes $j1,\ldots ,jL_j 
)$. Specifically, we use a linear map from $R^{\vert L\vert }$ to $R^{\vert 
L\vert -J}$ represented by a block-diagonal matrix $A$ with $J$ blocks of 
size $L_j \times (L_j -1)$:
    \begin{equation}
    \label{eq16}
        A_j = \left(
	\begin{array}{l}
	\begin{array}{ccccc}
                - \frac{\sqrt{L_j}-1}{L_j-1}   & \;\;\;1\;\;\; & \;\;\;0\;\;\; 
		& \dots & \;\;\;0 
	\end{array}
	\\
	 \cdots \cdots \cdots \cdots \cdots \cdots \cdots \cdots \cdots 
	 \\ 
	\begin{array}{ccccc}
                - \frac{\sqrt{L_j}-1}{L_j-1}   & \;\;\;1\;\;\; & \;\;\;0\;\;\; 
		& \dots & \;\;\;0 
	\end{array}
              \end{array}
	      \right)
    \end{equation}
Geometrically, such a map provides isometric rotation ($\bar {\lambda 
}^k=A\lambda ^k)$ to the hyperplane with zero first coordinate, i.e., (every 
block $A_j $ defines a rotation of a unit simplex in $L_j $--dimensional 
space around a hypersurface opposite to the first vertex; the angle of this 
rotation is such that the first vertex moves to the point where the first 
coordinate equals 0). Explicitly, this rotation is $\bar {\lambda }_{jl-1}^k 
=A_j \lambda _{jl}^k $ in matrix form or $\bar {\lambda }_{jl-1}^k =\lambda 
_{jl}^k -\frac{\sqrt {L_j } -1}{L_j -1}\lambda _{j1}^k $ for $l=2,\ldots 
,L_j $. New vectors $\bar {\lambda }^k $ do not possess any ties. It is easy 
to ascertain that such a transformation really conserves distances between 
points in a simplex. The reverse transformation is, 
\begin{equation}
\label{eq17}
\lambda _{j1}^k =\frac{1-\sum\nolimits_{l=2}^{L_j } {\bar {\lambda 
}_{jl-1}^k } }{\sqrt {L_j } },\,\,\,\,\,\,\lambda _{jl}^k \,=\bar {\lambda 
}_{jl-1}^k +\frac{\sqrt {L_j } -1}{L_j -1}\lambda _{j1}^k .
\end{equation}

\subsubsection{Algorithm for identifying the subspace}
The initial completion of 
the moment matrix is constructed in a arbitrary way, e. g, by the unitary 
diagonal matrix or completing by frequencies as $f_{ij} =f_i f_j $. The next 
preliminary step is the rotation of each simplex (corresponding to each 
question as described above) to the hyperplane to eliminate one degree of 
freedom.  This produces $n$ 
points $c^1 ,\ldots ,c^n $(images of columns of frequency 
matrix) in $m=(\vert L\vert -J)$­dimensional space. The problem is to find an 
affine plane that minimally deviates from these points in the space of 
individual probabilities. First, we find the center of gravity of this 
system 
\begin{equation}
\label{eq18}
c^0 =\frac{1}{n}\sum\limits_i {c^i } ,
\end{equation}
and then consider a new set of points $\bar {c}^i=c^i-c^0$, that corresponds 
to shifting the point of origin. Now we need to find a $K$-dimensional 
linear subspace in $R^m$ that minimally deviates from this set of points. 
The solution of this problem is well-known (see, e.g., chapter 43 of 
Kendall and Stuart, 1977): one has to consider an $m\times m$ matrix $X$ 
with components$X_{rs} =\sum\nolimits_i {\bar {c}_r^i } \bar {c}_s^i $; this 
matrix is symmetric and positively defined, and thus its normalized 
eigenvectors are composed of an orthonormal basis in $R^m$. Let $\gamma _1 \ge 
\gamma _2 \ge \cdots \ge \gamma _m >0$ be eigenvalues of matrix $X$, and let 
$z^1,\ldots ,z^m$ be corresponding eigenvectors. The plane of 
dimensionality $K$ that minimizes the sum of squared distances from points 
$\bar {c}^1,\ldots ,\bar {c}^n$ is spanned by $z^1,\ldots ,z^m$, and the sum 
of squared distances is $\mbox{tr }X-\sum\nolimits_{k=1}^K {\gamma _k } $. 
Vectors $c^0,c^0+z^1,\ldots ,c^0+z^{K-1}$ give us an affine basis of the 
sought affine plane. Finally, we apply inverses of transformation (\ref{eq17}) 
to $c^0,c^0+z^1,\ldots 
,c^0+z^{K-1}$ to obtain the sought basis $\lambda ^1,\ldots ,\lambda ^K$ of 
the subspace $\Lambda $.

\subsection{Choice of a basis}

The basis cannot be defined uniquely, and any convex combination keeping the 
LLS restrictions can be considered an alternative. 
A choice may be made 
using prior information about the process of interest. The appeal of 
prior information at this stage is reasonable because of the evident fact that 
the same dataset can be used for analyzing different (say, disability or 
CVD) processes. 

The way how this information is used and how the procedure of specific choice of
the basis is defined is a question of taste. We describe here two possible
schemes implemented in the algorithm prototype.

A researcher specifies the characteristics 
of ``ideal'' individuals based on his/her experience in the research domain. 
Then he/she can construct vectors of probabilities for such ideal individuals or 
take these individuals from the sample under consideration. The vectors of 
probabilities of these individuals are taken as basis vectors. If 
probability vectors are constructed by hand, they could be beyond 
polyhedron $P_g$, so they should be projected to $P_g$. The 
individual coordinates in this basis would represent ``proximity" of the 
individual to the ``ideal" ones.

In another scheme, the basis is obtained using assignment of LLS 
scores (calculated on some arbitrary basis) to $K$ clusters, and then basis 
vectors $\lambda ^1,\ldots ,\lambda ^K$ are calculated as means of 
probabilities $\beta _{jl}^i $ over each cluster.

A researcher can develop his/her own scheme of basis selection. For example,
he/she can simply use vectors already known 
from previous studies or construct a basis purely mathematically, e.g., from 
the condition of maximal 
linear independence of the vectors, or choose it from the 
set of the supportive polyhedron vertexes. 

\subsection{Calculation of individual conditional expectations. }

When a basis of the supporting plane is found, the conditional expectations 
can be found from the main system of equations (\ref{eq10}), which is a linear 
system after substituting the basis. The system, however, relates conditional
expectations $\Expec(G_k \vert X=\ell )$ for a pattern $\ell$ with at least one
0 outcome. Thus exact system of equations (\ref{eq10}) can be written for all
patterns $\ell$ except patterns where all outcomes are known. For the
complete patterns, we can calculate $J$
conditional expectations, subsequently excluding one of $J$ questions (i.e.,
obtaining patterns $\ell^{[j]}$, where $\ell^{[j]}$ denotes 
vector $\ell$ with $j^\text{th}$
coordinate equal to 0), solving
the exact system of equations for obtained patterns, and defining LLS score for
complete pattern as mean over $J$ solutions for conditional expectations 
for $\ell^{[j]}$ patterns. 
This approach can be formalized by considering a system of 
$J$ system of equations:
\begin{equation}
\label{eq19}
\sum\limits_k {\lambda _{jl}^k } \cdot g_{\ell k} \approx \frac{f_{\ell} }{f_\ell^{[j]} }
\end{equation}
 This is a sparse 
overdetermined system; methods for solving such systems are well-elaborated 
(see, for example, Forsythe et al. (1977), Kahaner et al. (1988)). Since 
$\sum\nolimits_k {g_{\ell k} } =1$ and $\sum\nolimits_k {\lambda _{jl}^k } 
\cdot g_{\ell k} \ge 0$ have to be satisfied for any $\ell $, the 
Lagrangian for the optimization task has to include corresponding terms with 
Lagrange multipliers:
\begin{equation}
\label{eq20}
Lg_\ell =\sum\limits_j {\left( {\sum\limits_k {\lambda _{jl}^k } \cdot 
g_{\ell k} -\frac{f_{\ell } }{f_\ell^{[j]} }} \right)} ^2+\sum\limits_{jl} 
{\rho _{jl} \sum\limits_k {\lambda _{jl}^k } \cdot g_{\ell k} } +\rho \left( 
{\sum\limits_k {g_{\ell k} } -1} \right)
\end{equation}
This is a quadratic programming problem. We used the method of 
antigradient projection to the 
polyhedron defined by the inequality and equality conditions with iterations 
as described in Section 3.7. The initial point is defined as $g_{\ell k} =1/K$. 
It is possible to show that this initial point always belongs to the 
polyhedron. 

LLS scores can deviate restrictions $\sum\nolimits_k {\lambda _{jl}^k } \cdot 
g_{\ell k} \ge 0$ in limited numbers of cases; e.g., for a relatively 
small number of questions $J<100$. Options to solve the optimization task 
without the inequality restriction are also included in the program. 
Keeping in (\ref{eq21}) Lagrange multipliers for equality conditions only we have 
\begin{equation}
\label{eq21}
Lg_\ell =\sum\limits_j {\left( {\sum\limits_k {\lambda _{jl}^k } \cdot 
g_{\ell k} -\frac{f_{\ell } }{f_{\ell^{[j]}} }} \right)} ^2+\rho \left( 
{\sum\limits_k {g_{\ell k} } -1} \right)
\end{equation}
The advantage of this task is that it can be solved explicitly, i.e.; without the 
search procedures required for the quadratic programming problem. The 
simplest way to deal with the restriction is to just substitute the $K^{th}$ 
component and apply SVD to a $J\times (K-1)$-matrix $\Lambda _{jk} =\lambda 
_{j\ell _j }^k -\lambda _{j\ell _j }^K $. If SVD has explicit form as 
$\Lambda _{jk} =\sum\nolimits_{k'} {\bar {\Lambda }_{jk'} } W_{k'} V_{k'k} 
$, where $\bar {\Lambda }$ is $J\times (K-1)$-matrix, $W$ contains a vector 
of singular values, and $V$ is an orthogonal matrix, then the result for 
conditional expectations is
\begin{equation}
\label{eq22}
g_{\ell k} =\sum\limits_{k'=1}^{K-1} {\frac{V_{kk'} }{W_{k'} }} \bar 
{\Lambda }_{jk'} \left( {\raise0.7ex\hbox{${f_{\ell} }$} 
\!\mathord{\left/ {\vphantom {{f_{\ell +l_j } } {f_{\ell^{[j]}} 
}}}\right.\kern-\nulldelimiterspace}\!\lower0.7ex\hbox{${f_{\ell^{[j]}} }$}-\lambda 
_{jl}^K } \right)_j 
\end{equation}
for $k=1,\ldots ,K-1$ and $g_{\ell K} =1-\sum\nolimits_{k=1}^{K-1} {g_{\ell 
k} } $.

\subsection{Mixing distribution}

The mixing distribution for an analyzed set of data is approximated 
by empirical distribution, where an individual gives a unit contribution to 
the histogram of the distribution. A support of this distribution is 
a set of $I$ points. Probabilities of the joint distribution (\ref{eq4}) 
are estimated as the sum over sample individuals or to the sum over possible
outcome patterns,
\begin{equation}
p_\ell^\ast 
=\sum\nolimits_i {\prod\nolimits_{j:\ell _j \ne 0} {\beta _{j\ell _j }^i } } 
=\sum\nolimits_{\ell '} f_{\ell '} 
\prod\nolimits_{j:\ell _j \ne 0} \sum\nolimits_k g_{\ell 'k} \lambda _{j\ell _j
}^k  
.
\end{equation}

\subsection{Quadratic programming problem.}

The quadratic programming problem is used by the algorithm in the completion 
procedure, in calculation of individual LLS scores, and in basis construction. 
Corresponding expressions for Lagrangians are given in eqs. (\ref{eq15}), (\ref{eq20}), and 
(\ref{eq23}). In all cases, the function minimized is quadratic, and equality and 
inequality conditions are linear functions of the sought parameters. The matrix of 
quadratic form is positive definite in all cases. The solution is based on 
Kuhn-Tucker conditions which become sufficient in these cases. The method 
of antigradient projection is used in searching for the optimal point (Attetkov et 
al., 2001). The antigradient direction is defined by the direction opposite 
to the gradient in a current point. Then the antigradient is projected to a subspace 
defined by so-called active conditions: all equality conditions and those 
inequality conditions which are strictly satisfied (i.e., the current point 
belongs to the corresponding plane). If the absolute value of the projection is not 
zero then the current point is moved in this direction up to the 
intersection with the closest boundary, or to the minimum of the 
quadratic function, whichever is closer. If the absolute value of the projection 
is zero, then the sign of Lagrange multipliers corresponding to active 
inequality restrictions is checked. If they all are positive, the 
current point is the minimum; otherwise a restriction corresponding 
to the negative Lagrange multiplier is removed and a 
new antigradient projection 
 is calculated.

\subsection{Clustering}

Clustering LLS scores is not a necessary component for  finding the 
subspace or calculation of LLS scores; however, it is helpful in selecting a
basis and in cross-checking in simulation studies. Two types of clustering 
procedures (Manton et al., 2003) are implemented in the algorithm. The 
first is the k-means algorithm for the 
situation where the number of clusters is fixed a priori. The second 
is a hierarchical procedure, which 
sequentially joints clusters with minimal distance between them. Several 
distance definitions are possible: distance between centers of mass 
in clusters, between closest and the most outlying cluster members. 
Numerical comparison shows that the most reliable results are obtained for 
the center-of-mass scheme. 

\subsection{Missing data}\label{sectmissing}
Missing data are often difficult problems in statistical 
analyses. Missing data are generated by the absence of responses for an 
individual to specific questions. The properties of LLS analysis make this 
kind of missing data relatively easy to handle. Two main sources of missing 
answers could be considered: first, when the failure to answer the question 
is random; and second, when the failure to answer the question correlates 
with answers to other questions. In the first case (missing data are 
random), the solution is provided by the fact that the input of the LLS 
algorithm consists of frequencies of partial response patterns (like the 
frequency of giving response C to the 2$^{nd}$ question and response A to 
the 5$^{th}$ question). With missing data, such frequencies can be 
calculated by relating the number of individuals with a particular 
response pattern to the number of individuals who gave answers to the 
questions covered by the response pattern (rather than to the total number 
of individuals). The only drawback of this method is decreased precision 
of frequency estimators. As LLS scores are expectations of latent variables 
conditional on the arbitrary part of the response pattern for an individual, 
the available part of the response pattern can be used to estimate the value 
of the latent variable. In the second case (missing data are not random), 
the absence of an answer can be considered an additional alternative for 
answering a multi-choice question; in this case, standard LLS analysis can be 
applied.

\section{Applications}

Below we i) perform a series of simulation experiments to test the 
predictive power of LLS model and its ability to reveal and to 
quantitatively reconstruct a hidden latent structure and ii) demonstrate how the
model performs when applied to a real dataset considering as an example a
demographic survey with 57 questions about disability and health status in 
5,161 individuals (Manton and Gu, 2001).

\subsection{Simulation studies}
 We perform three types of simulation experiments which 
demonstrate the quality of the reconstruction of i) linear subspace, ii) LLS 
score distribution, and iii) clusters in LLS score space. At the first stage 
of each simulation study, the dimensionality of the latent structure 
(strictly, dimensionality of the linear space supporting the structure in 
the space of individual probabilities) is chosen. Then, assuming a specific 
structure of questions and numbers of outcomes to each question (i.e., set 
of indexes$j\ell _j )$, basis vectors $\lambda ^1,\ldots ,\lambda ^K$and 
mixing distribution $G$ are constructed from some arguments or simulated 
randomly. Then using individual probabilities calculated as in Eq. (\ref{eq2}), 
outcomes of the specified number of individuals are simulated. The 
database generated is used as an input for LLS analysis to extract the supporting linear 
$K$-dimensional subspace and to estimate the mixing distribution via estimation 
of individual conditional expectations. The comparison of original and 
extracted linear subspace and mixing distribution provides information about 
the quality of LLS analysis in  conditions defined by $K$, $I$, $J$, 
and $L_j $.

\subsubsection{Linear subspace} 

We started with investigation of the quality of reconstruction of 
linear subspace. Choosing initial 2-, 3-, and 5-dimensional subspaces with 
60, 120, and 240 dichotomous questions and assuming uniform distribution of LLS 
scores, we simulated 1,430 and 14,300 sets of individual responses. Then we 
reconstructed a linear subspace and compared it with the initial linear subspace. For 
such comparisons, we needed a procedure for the quantitative comparison of 
closeness of multidimensional subspaces, i.e, a measure of a distance or 
angle between the linear subspaces. Such measure $d$ is constructed on the 
basis of the theory of orthogonal projector operators (see Golub and Van 
Loan, 1989), by which it is possible to define a distance between linear 
subspaces. By using CS-decomposition, it was proved that the distance may be 
interpreted as the sine of the generalized angle between subspaces. 
Technically, $d=\sqrt {1-\sigma _{\min }^2 } $, where $\sigma _{\min } $ is the 
minimal singular value obtained by SVD of $K$-dimensional matrix $Q=P_1^T 
P_2 $, and $P_1 $ and $P_2 $ are $\vert L\vert \times K$- orthogonal 
projectors onto comparing subspaces. We will use $d$ as a measure of the quality 
reconstruction.  

Specifically, original vectors $\lambda ^1,\ldots ,\lambda ^K$ are 
constructed as follows: i) vectors $\lambda ^1,\ldots ,\lambda ^K$ can only 
take values 0 or 1, ii) $\lambda _{j1}^1 =1$ for all $j$, iii) for $k\ge 2$ 
$\lambda _{j1}^k =0$ if $j\le J \mathord{\left/ {\vphantom {J {(K-1)}}} 
\right. \kern-\nulldelimiterspace} {(K-1)}$ and $\lambda _{j1}^k =1$ 
otherwise. In all cases $\lambda _{j2}^k =1-\lambda _{j1}^k $. The distribution 
of $G$ was chosen to the concentrated in a set of discrete points uniformly covering 
the $K$-dimensional simplex. Then a set of probabilities $\beta _{jl}^i $ is 
calculated for all individuals and $J$ outcomes per individual are simulated 
according to these probabilities. To get stable results we 
performed a sufficient number of experiments (on the order of 100) and obtained 
$d$ presented in the Table \ref{tab1}.

\begin{table}[htbp]
\caption{Distances between simulated and reconstructed subspaces}
\begin{tabular}{ccccccccc}
\hline
& 
\hspace{0.25in} 
&
\multicolumn{3}{c}{$N$=1,430} 
&
\hspace{0.25in} 
&
\multicolumn{3}{c}{$N$=14,300} 
 \\
\cline{3-5}
\cline{7-9}
$K$&& 
$J$=60 & 
$J$=120 & 
$J$=240 & 
&
$J$=60 & 
$J$=120 & 
$J$=240  \\
\hline
2&& 
0.023& 
0.023& 
0.022& 
&
0.008& 
0.007& 
0.007 \\
3&& 
0.075& 
0.072& 
0.070& 
&
0.023& 
0.023& 
0.023 \\
5&& 
0.222& 
0.190& 
0.176& 
&
0.073& 
0.059& 
0.057 \\
\hline
\end{tabular}
\label{tab1}
\end{table}

\subsubsection{LLS score distribution}

We designed the experiments to have statistical stability and to 
analyze a large number of questions ($J=1,500)$. The position of the 
supporting plane was chosen randomly. LLS scores are distributed over a union 
of intervals $[0.10,0.25]$ and $[0.50,0.75]$. Specific values of LLS scores 
are chosen deterministically but not simulated as in Kovtun et al. (2005b): 
5,000 from 0.10 to 0.25 and 5,000 from 0.50 to 0.75. In total, a sample of 
10,000 individuals was generated. First, a mixing distribution was restored 
from the sample data using the original basis used for outcome simulation. 
Figure 1a demonstrates reconstructed and true distributions. The latter is 
shown as a solid line. The figure demonstrates an excellent 
quality of restoration of the mixing distribution when for reconstruction of 
mixing distribution we used a simulated, but not reconstructed, basis. Using 
original vectors $\lambda ^1$ and $\lambda ^2$ in this comparison allows us 
to test the procedure of reconstruction of LLS scores separately. These 
results together with results of a separate test of quality of supporting 
subspace reconstruction in Table 1 complete the test of the two 
main components of the LLS algorithm.

However, it is still not clear how overall reconstruction is good, i.e. how 
uncertainties of the subspace reconstruction could impact the reconstruction of
the 
mixing distribution. To perform this test we used reconstructed vectors 
$\lambda ^1$ and $\lambda ^2$ instead of original in the example. The 
results are given in Figure 1b.

\begin{figure}[htb]
\includegraphics[height=7.2cm]{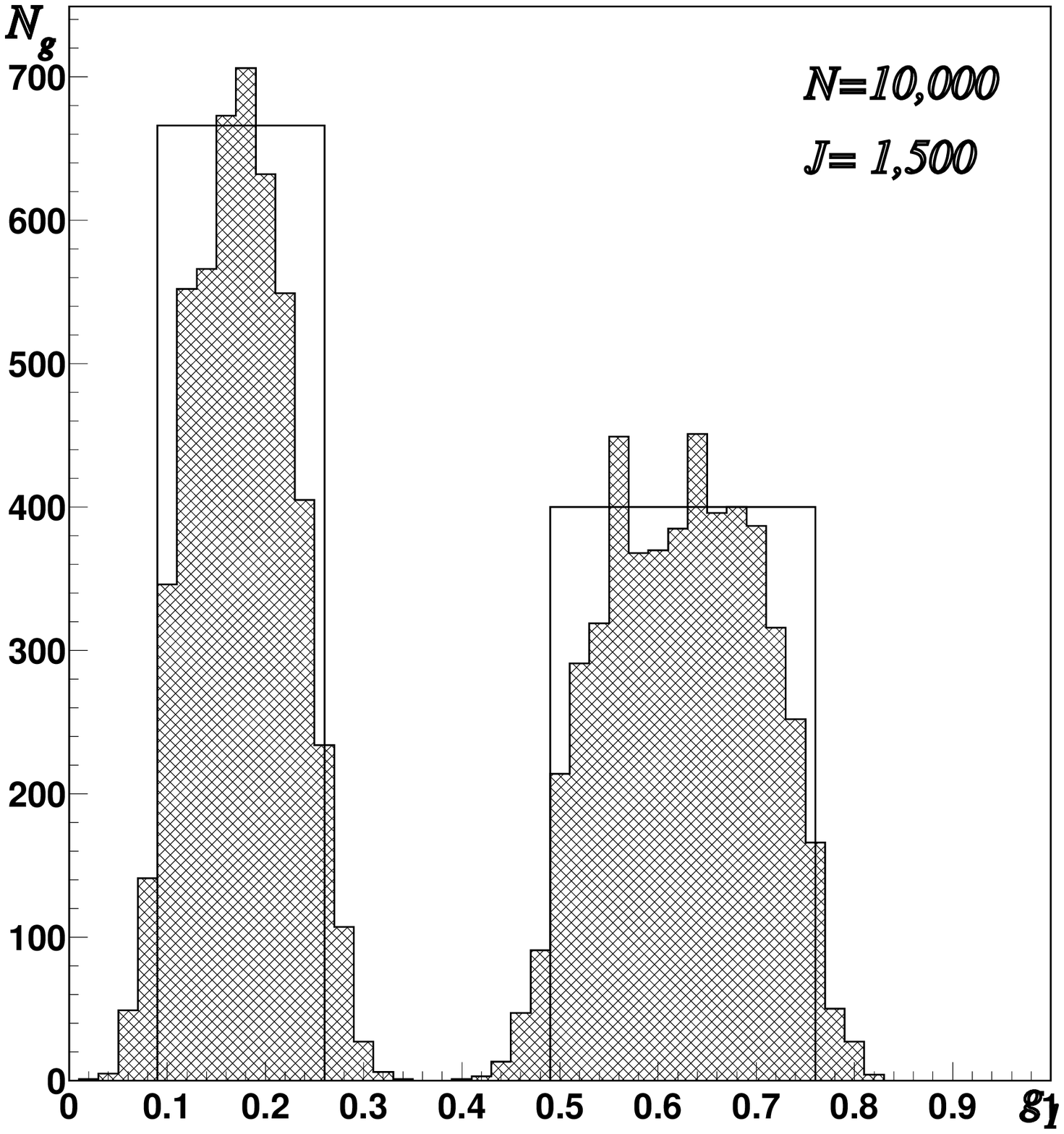}
\includegraphics[height=7.2cm]{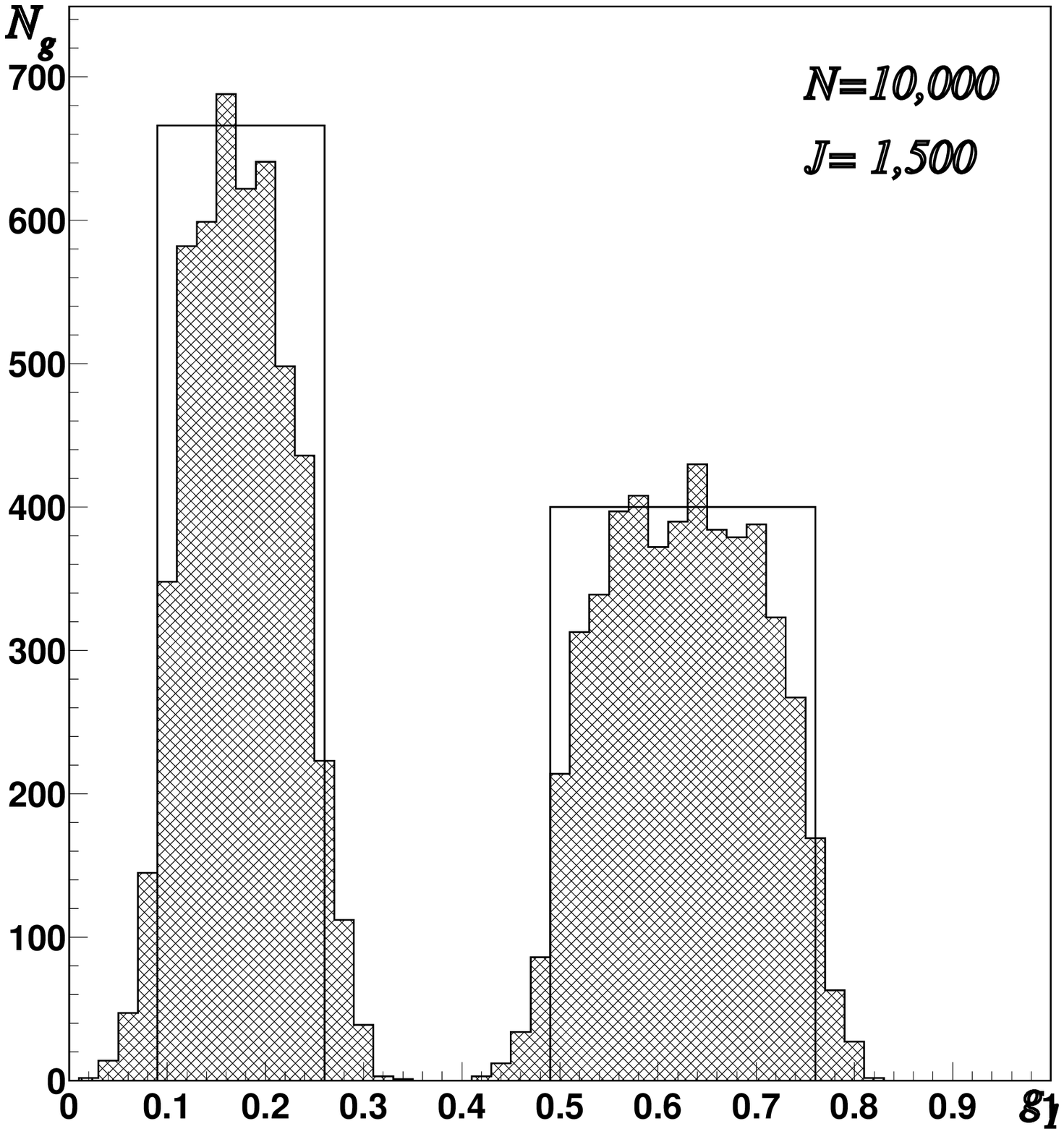}
\caption{Reconstructed mixing distribution using simulated (left) and
reconstructed (right) basis. Lines show 
simulated mixing distribution}\label{fig1}
\end{figure}

In addition, to analyze the performance of the algorithm we conducted this 
experiment for different numbers of individuals $I$ and variables $J$. The 
time spent on the algorithm almost does not depend on the size of the 
sample, as the most time-consuming part of the algorithm is SVD of moment 
matrix, which has to be performed 6 times (for rank estimation and for 5 
iterations of supporting subspace finding). Size of moment matrix depends on 
the number of measurements $J$. In the table below we give the time spent by 
the algorithm for different values of $J$ using a computer with an AthlonXP 
2500 (2.5GHz) processor and 512MB DDR RAM. In all cases sample size is 8100.

\begin{table}[htbp]
\caption{Running time for computer experiments with different $J$}
\begin{tabular}{*{10}c}
\hline
$J$&\hspace{0.1in}& 
50& 
100& 
200& 
300& 
500& 
750& 
1000& 
1500 \\
\hline
Time&& 
18s& 
23s& 
1m 49s& 
4m 45s& 
17m 40s& 
59m& 
2h 40m& 
10h 35m \\
\hline
\end{tabular}
\label{tab2}
\end{table}

\subsubsection{Clusters in LLS score space} We checked the ability of the LLS 
analysis to identify clusters of LLS scores. We assumed the LLS score 
distribution was concentrated at five points with equal weights. Then 
setting $K=3$ and $I=1,000$ we simulated individual responses for $J$=100, 
200 and 500. Then we applied a complete-linkage hierarchical clustering 
procedure 
(Section 3.8) both to the reconstructed LLS scores and to the original vector of 
response to combine individuals in 5 classes. As a measure of quality, we 
used the fraction of individuals assigned to the correct class. We performed 
10 simulations to make the calculation statistically stable. The average 
number of misclassified individuals is given in Table \ref{tab3}. 

\begin{table}[htbp]
\caption{The average number of misclassified individuals for different $J$}
\begin{tabular}{p{110pt}rrr}
\hline
$J$& 
100& 
200& 
500 \\
\hline
Original responses& 
22.7{\%}& 
17.3{\%}& 
14.9{\%} \\
LLS scores& 
4.8{\%}& 
0.2{\%}& 
0.0{\%} \\
\hline
\end{tabular}
\label{tab3}
\end{table}

This example demonstrates that the individual conditional expectations are 
better suited for the purpose of classification then the original vectors of 
responses. 

\subsection{Application to the NLTCS data}

The National Long Term Care Survey is a longitudinal survey designed to 
study changes in the health and functional status of older Americans (aged 
65+). It also tracks health expenditures, Medicare service use, and the 
availability of personal, family, and community resources for caregiving. 
The survey began in 1982, with follow-up surveys done in 1982, 1984, 
1989, 1994, and 1999. A sixth follow-up survey was conducted in 2004-2005. A 
detailed description of the NLTCS may be found at \underline 
{http://nltcs.cds.duke.edu/}.

We analyzed a sample of 5,161 individuals from the 1999 NLTCS, and selected 
57 questions. 27 variables characterize disability level 
with respect to activities of daily living, instrumental activities of daily 
living, and physical impairment. 30 variables describe self-reports 
of chronic diseases. Details about these questions may be found in 
Manton and Gu (2001). Then we excluded 370 individuals with at 
least one missing outcome. Thus, 4,791 remaining individuals are subject to 
57 questions, 49 of which have 2 possible answers, and 8 have 4 possible 
answers. In total, we have $\left| L \right|=130$.

The first task is to define the dimensionality of the LLS problem, $K$. As 
 mentioned above, $K$ is the rank of the moment matrix, and thus it 
could be estimated as the rank of the frequency matrix. We used the SVD to 
estimate the rank of the frequency matrix. As the frequency matrix is 
incomplete, we decompose not the whole matrix, but 
its bottom left corner of size 64. The singular values obtained were compared 
to total statistical error found as a quadratic sum of cell errors. This ratio 
was 0.292. Individual cell standard errors were estimated using 
Wilson's approach (\ref{eq14}). The first 10 singular values are given in Table
\ref{tab4}.

The computations show that the hypothesis $K=4$ has to be accepted with 
confidence level 95{\%}; that corresponds approximately to the double standard 
error interval (0.584). For comparison, we will also consider the case of 
$K=3$. Note that there is no significant gap between the 4th and 5th singular 
numbers. This suggests that the support of the distribution is an ellipsoid of 
full dimensionality which is thinner in higher dimensions, and choosing a 
particular value for $K $ approximates this ellipsoid by a lower-dimensional 
ellipsoid obtained from the true one by collapsing a number of smaller axis.

\begin{table}[htbp]
\caption{First 10 singular values of frequency matrix of NLTCS}
\begin{tabular}{*{15}c}
\hline
$\sigma _1 $& 
$\sigma _2 $& 
$\sigma _3 $& 
$\sigma _4 $& 
$\sigma _5 $& 
$\sigma _6 $& 
$\sigma _7 $& 
$\sigma _8 $& 
$\sigma _9 $& 
$\sigma _{10} $
\\
\hline
39.112 & 
3.217& 
1.464& 
0.652& 
0.363& 
0.310& 
0.243& 
0.220& 
0.198& 
0.148
\\
\hline
\end{tabular}
\label{tab4}
\end{table}

When the dimensionality of the LLS-problem is fixed, we can complete the moment 
matrix using the algorithm described in Section 3.3. The sub-matrix 
corresponding to 
the first four dichotomous variables is,
\[
\left( {\begin{array}{l}
 \mbox{ 0.094 }\,\,{\rm {\bf 0.513 }}\,\,{\rm {\bf 0.051}}\mbox{ 
}\,\,\mbox{0.328 }\,\,\mbox{0.011 }\,\,\mbox{0.258 }\,\,\mbox{0.012 
}\,\,\mbox{0.518 }\,\,\mbox{0.014 } \\ 
 \mbox{ 0.906 }\,\,{\rm {\bf 0.487 }}\,\,{\rm {\bf 0.949}}\mbox{ 
}\,\,\mbox{0.672 }\,\,\mbox{0.989 }\,\,\mbox{0.742 }\,\,\mbox{0.988 
}\,\,\mbox{0.482 }\,\,\mbox{0.986 } \\ 
 \mbox{ 0.264 }\,\,\mbox{0.918 }\,\,\mbox{0.196}\,\,\mbox{ }{\rm {\bf 0.633 
}}\,\,{\rm {\bf 0.128}}\mbox{ }\,\,\mbox{0.688 }\,\,\mbox{0.051 
}\,\,\mbox{0.846 }\,\,\mbox{0.153 } \\ 
 \mbox{ 0.736 }\,\,\mbox{0.082 }\,\,\mbox{0.804 }\,\,{\rm {\bf 
0.367}}\,\,{\rm {\bf 0.872}}\mbox{ }\,\,\mbox{0.312}\,\,\mbox{ 0.949 
}\,\,\mbox{0.154}\,\,\mbox{ 0.847 } \\ 
 \mbox{ 0.335 }\,\,\mbox{0.916 }\,\,\mbox{0.275 }\,\,\mbox{0.872}\,\,\mbox{ 
0.142 }\,\,{\rm {\bf 0.664}}\,\,{\rm {\bf 0.164}}\,\,\mbox{ 
0.888}\,\,\,\,\mbox{0.230 } \\ 
 \mbox{ 0.665}\,\,\mbox{ 0.084 }\,\,\mbox{0.725 }\,\,\mbox{0.128}\,\,\mbox{ 
0.858 }\,\,{\rm {\bf 0.336 }}\,\,{\rm {\bf 0.836}}\mbox{ 
}\,\,\mbox{0.112}\,\,\mbox{ 0.770 } \\ 
 \mbox{ 0.160 }\,\,\mbox{0.879 }\,\,\mbox{0.085 }\,\,\mbox{0.514 
}\,\,\mbox{0.034 }\,\,\mbox{0.424 }\,\,\mbox{0.027 }\,\,{\rm {\bf 0.640 
}}\,\,{\rm {\bf 0.069}}\mbox{ } \\ 
 \mbox{ 0.840 }\,\,\mbox{0.121 }\,\,\mbox{0.915 }\,\,\mbox{0.486 
}\,\,\mbox{0.966 }\,\,\mbox{0.576 }\,\,\mbox{0.973 }\,\,{\rm {\bf 0.360 
}}\,\,{\rm {\bf 0.931}}\mbox{ } \\ 
 \end{array}} \right)
\]
Completed values are marked in bold style. Renormalization of the frequency 
matrix is performed (see Section 3.1). Recall that columns of the frequency 
matrix are approximately in the sought linear subspace. The procedure 
described in Section 3.3 provides an unambiguous definition of the sought 
subspace and some basis vectors $\lambda _{jl}^k $. The basis cannot be 
defined uniquely, and any convex combination keeping LLS restrictions (3) 
can be considered an alternative. In this example, we apply a combination of 
the methods discussed in Section 3.4. Specifically, in the 
beginning we estimated LLS scores in some prior basis. Since subsequent 
conclusions are made on the basis of clusterization and analysis of 
probability vectors, which are independent of initial basis choice, the 
initial basis is taken as it was obtained from subspace construction, i.e., 
it is arbitrary. Then LLS scores are assigned to 7-8 clusters (i.e., for slightly 
larger numbers of clusters than $K)$ using methods of cluster analysis. By 
analyzing outcomes of typical representative respondents (with LLS scores 
closest to centers of clusters) and of probabilities $\beta _{jl}^i $ 
calculated as means over each cluster, we can identify the cluster structure of 
the analyzed population. Based on this analysis, we chose $K=4$ clusters 
corresponding to i) healthy individuals ($k=1)$, ii) partly disabled 
individuals ($k=4)$, iii) strongly disabled individuals ($k=2)$, and iv) 
individuals with chronic diseases but without evidence of disability 
($k=3)$. For each such group, we created a ``typical'' vector of 
probabilities $\bar {\beta }_{jl}^i $ by hand. Specifically, for the first 
group we chose unit probabilities for all answers corresponding to healthy 
states. For the second group, mean frequencies over the sample are chosen. 
For the third group, unit probabilities are assigned to answers 
corresponding to strongly disabled states and means over samples for chronic 
disease questions. For the fourth group, unit probabilities correspond to 
non-disabled states and to all disease diagnoses. The final basis is constructed 
by projection of all four vectors to the simplex. Technically, finding 
vectors $g_k $ minimizes the Lagrangian (\ref{eq20}), which for this task has the 
form:
\begin{equation}
\label{eq23}
Lg_{k'} =\sum\limits_{jl} {\left( {\sum\limits_k {\lambda _{jl}^k } \cdot 
g_k -\beta _{jl}^{k'} } \right)} ^2+\sum\limits_{jl} {\rho _{jl} 
\sum\limits_k {\lambda _{jl}^k } \cdot g_k } +\rho \left( {\sum\limits_k 
{g_k } -1} \right),
\end{equation}
and applying the algorithm sketched in Section 3.7. Then a similar analysis 
was performed for $K=3$, where the basis component for partly disabled 
individuals was rejected. Figure 2 demonstrates the results. The plots on the left 
show 2d-polyhedron for$K=3$ and two projections of 3d-polyhedron for $K=4$. 
These polyhedrons are defined by LLS restrictions (3). In this case, LLS 
scores are restricted by 130 inequality and 1 equality. 2d-polyhedron is 
shown as it is in the plane normal to vector 
$g=(\raise0.5ex\hbox{$\scriptstyle 
1$}\kern-0.1em/\kern-0.15em\lower0.25ex\hbox{$\scriptstyle 
3$},\raise0.5ex\hbox{$\scriptstyle 
1$}\kern-0.1em/\kern-0.15em\lower0.25ex\hbox{$\scriptstyle 
3$},\raise0.5ex\hbox{$\scriptstyle 
1$}\kern-0.1em/\kern-0.15em\lower0.25ex\hbox{$\scriptstyle 3$})$. 
3d-polyhedron is in the plane normal to $g=(\raise0.5ex\hbox{$\scriptstyle 
1$}\kern-0.1em/\kern-0.15em\lower0.25ex\hbox{$\scriptstyle 
4$},\raise0.5ex\hbox{$\scriptstyle 
1$}\kern-0.1em/\kern-0.15em\lower0.25ex\hbox{$\scriptstyle 
4$},\raise0.5ex\hbox{$\scriptstyle 
1$}\kern-0.1em/\kern-0.15em\lower0.25ex\hbox{$\scriptstyle 
4$},\raise0.5ex\hbox{$\scriptstyle 
1$}\kern-0.1em/\kern-0.15em\lower0.25ex\hbox{$\scriptstyle 4$})$. Basis 
vectors produced unit simplexes are labeled by numbers. Plots on the right 
demonstrate how the polyhedrons are filled by the considered population. 
For the 
filling we performed clustering to 1,000 clusters. Each point in the 
plots represents one cluster. The area of each point is proportional to the 
number of individuals assigned to the corresponding cluster. The exception is 
the point marked by open circles with a closed point inside. About half of 
the total population was assigned to this cluster. We restrict consideration 
to this illustration.

\begin{figure}[btp]
\includegraphics[height=19cm]{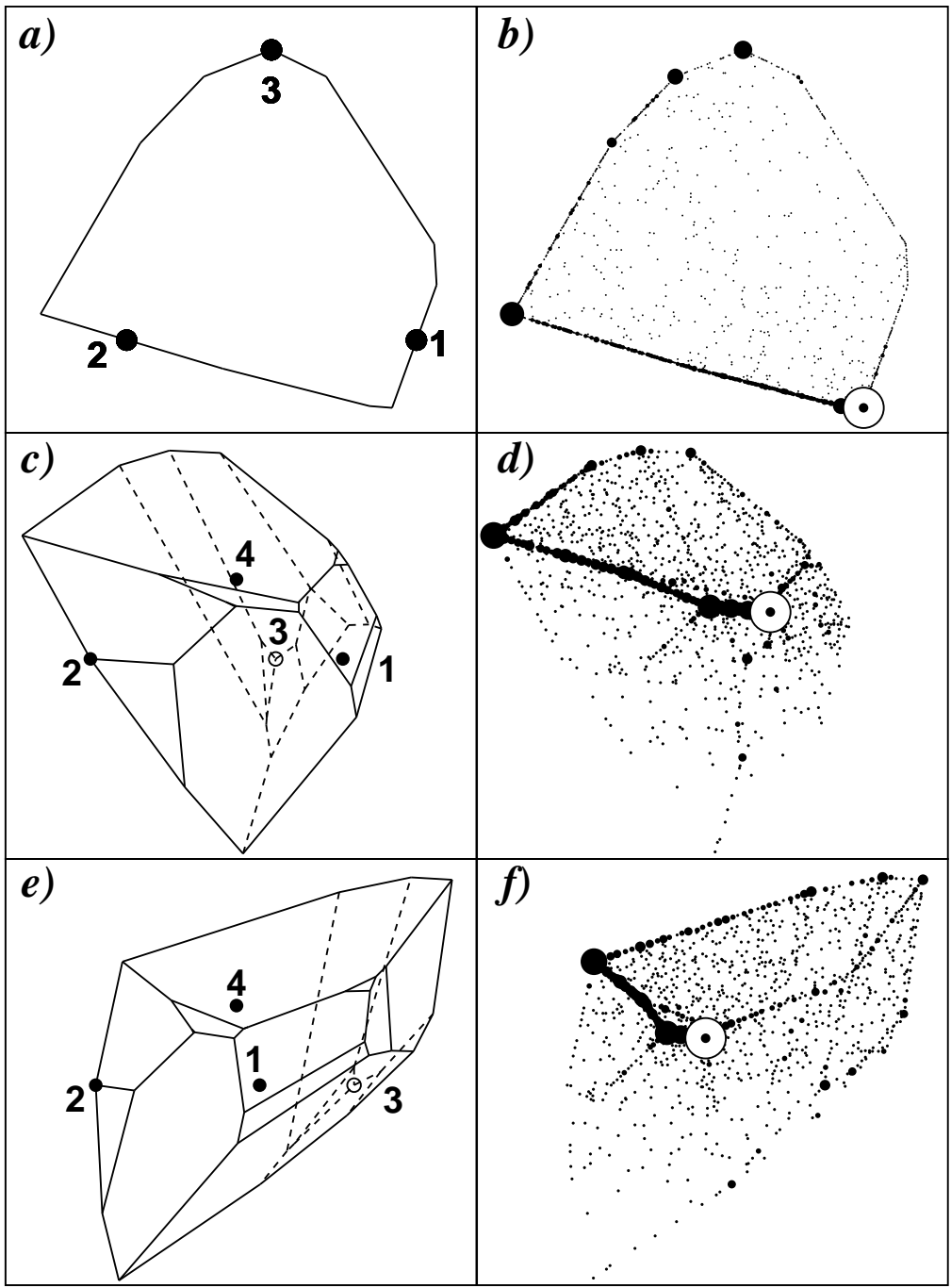}
\caption{Polyhedrons defined by LLS constrains for $K$=3 (a) and $K$=4 (c,e) and
their filling by LLS scores of NLTCS individuals (b,d,f). See text for further
explanations.
}\label{fig2}
\vspace{1cm}
\end{figure}


\section{Discussion and Conclusion}

LLS is a model describing high-dimensional categorical data assuming existence 
of a latent structure described by 
$K$-dimensional random vectors. This vector is interpreted as explanatory 
variables which can shed light on mutual correlations observed in measured 
categorical variables. This vector plays the role of a random variable mixing 
independent distribution such that the observed joint distribution is maximally close 
to the data. Mathematically, LLS analysis considers the observed 
joint distribution of categorical variables as a mixture of individual joint 
distributions, which are assumed to be  independent. This is the standard 
``local independence'' assumption of latent structure analysis; the 
specific LLS assumption is that the mixing distribution is supported by a 
$K$-dimensional subspace $\Lambda $ of the space of all independent 
distributions or, equivalently, of the space of individual probabilities. The 
mixing distribution can be considered as a distribution of random vectors 
$G$, which take values in $\Lambda $. The vectors $g_i $ (LLS scores) are the 
hidden states of individuals in which we are 
interested. They can be estimated as conditional expectations of $G$, $g_\ell 
=\Expec(G\vert X_1 =\ell _1 ,\ldots ,X_J =\ell _J )$, where $\ell =(\ell _1 
,\ldots ,\ell _J )$ is a response pattern. Support of this random vector is a 
$K$-dimensional space, the dimension of which is defined by the dataset. Linear 
support is a distinctive feature of LLS compared to other latent 
structure analyses. For example, LCM is characterized by mixing distribution 
concentrated in several isolated points. 

Another important distinction is the existence of an algorithm capable of 
estimating a LLS model for large numbers of questions and individuals. 
When a basis $\lambda ^1,\ldots ,\lambda ^K$ of linear subspace supported 
mixing distribution is 
known, conditional expectations $g_\ell $ can be calculated by solving a 
linear system of equations. A basis $\lambda ^1,\ldots ,\lambda ^K$, in turn, can be 
identified by applying principal component analysis to the moment matrix. 
As the choice of a basis is not unique, one has 
to apply substantial knowledge derived from the applied domain to make the 
most appropriate choice. This algorithm, being a sequence of linear 
algebra methods, does not use maximum likelihood methods. This is an 
advantage of the method, because individual information is presented via 
nuisance parameters, which creates difficulties in the maximum 
likelihood approach. The LLS algorithm for parameter estimation is based on 
two theorems which became possible because of a linear assumption about the 
support of the mixing distribution. The first identifies properties of the 
moment matrix. The second  presents the main system of equations. Just 
the existence of the system of equations that describes parameters of the 
model is a significant advantage of the LLS method, as it allows us to avoid 
using maximum likelihood estimators, which may not be consistent 
in the presence of nuisance parameters. 

Construction of LLS analysis on the basis of the theory of mixed 
distributions provides clear and simple proof of basic statistical 
properties of the estimator. Regarding consistency, note that parameter 
estimations of the LLS model are solutions of a quasilinear system of 
equations, whose coefficients are the observed frequencies. When true 
moments are substituted in the system instead of frequencies, its solution 
represents the true values of parameters. The consistency of LLS estimators 
follows from the fact that i) solutions of the system continuously depend on 
its coefficients, and ii) frequencies are consistent estimators for the 
moments. Regarding identifiability, it was shown in Kovtun et al. (2005b) 
that the LLS model is identifiable if the moment matrix is ``sufficiently 
non-degenerated''. It was also shown in Kovtun et al. (2005b) that for 
almost all mixing distributions the moment matrix is non-degenerated. Thus, 
LLS models are almost surely identifiable. 

We performed a series of simulation experiments to analyze the quality 
of reconstruction of: i) linear subspace; ii)  LLS 
mixing distribution; and iii) clustering 
properties. The results demonstrated an acceptable quality of reconstruction. It 
is interesting to compare the results of the LLS model with the potential 
results of a latent class model. If mixing distribution is concentrated in 
several isolated points, LCM would restore the same picture as LLS, and 
thus, it could be an alternative for the LLS model. However, in cases like 
in Figure 2 LCM may only be used as an approximation, and it may be shown 
that LCM should involve approximately 1,000 latent classes 
(the number of classes in this case is approximately equal to number of different 
LLS scores obtained by the LLS model).

Analysis of demographic data on disability and chronic diseases demonstrated 
the potential of the LLS model to investigate the properties of populations with latent 
heterogeneity. Performing cluster analysis, we identified four basic 
population groups corresponding to healthy, disabled, strongly disabled 
individuals, and persons with chronic diseases. We found linear 
subspaces for $K=3$ and $K=4$. Empirical mixing distributions for these cases 
are presented as distributions on supportive polyhedrons (Figure 2) defined 
by LLS restrictions. 

This analysis can be naturally generalized to the case of longitudinal data 
using a binomial quadratic hazard model (Manton et al, 1992) or stochastic 
process model (Yashin and Manton, 1997). This can provide the possibility of 
introducing mortality into the model. The basic advantages of these models 
are that they provide a biologically justified U- or J- shaped hazard function 
versus covariates (Witteman et al., 1994). In our case, LLS scores would 
play the role of covariates. LLS scores are constrained to have a sum of one. 
This is appropriate for interpretation and applications. Constraints of LLS 
scores (e.g., $\sum\nolimits_k {g_k } =1)$ prevent the direct use of 
a stochastic process model for projection. One way to eliminate this 
restriction without losing the advantage in interpretation is to eliminate 
the component corresponding to the healthiest pure type (due to the freedom 
in basis choice, a basis with healthy state can always be chosen), thereby 
solving the problem of dealing with this restriction and, 
moreover, allowing us to specify a model for the hazard function. For the 
hazard we will use a quadratic hazard model, where the quadratic form 
models dependence on the components of state vectors. Since the health 
component of the LLS score vector is eliminated, it is reasonable (and 
technically appropriate) to assume that the minimum of mortality is located 
in the origin over the remaining components (i.e., mortality is minimal for 
the healthiest people). This allows us to drop the linear component of the 
quadratic hazard function. After estimating the parameters, a scheme for 
the projection calculation has to be developed. In the general case, possibly 
including nonlinear relations in dynamic equations, the approach can be based on a 
microsimulation procedure; i.e., a method for simulating trajectories for each 
person in the cohort (Akushevich et al., 2005). 

Furthermore, LLS opens its own possibilities for longitudinal analyses of 
repeated measurements using its inherent features. Supporting planes 
obtained for different waves may not coincide with each other. The 
first question to be analyzed is how different are supportive planes obtained 
from different waves and whether a plane obtained for combined dataset is a 
good fit or there is a trend in the plane movement over a state space. In 
the first case, a plane obtained for joint analysis would be taken for 
longitudinal analysis. Basis becomes time independent and chosen by taking 
specifics of the problem into account. If a trend of time movement of 
supportive plane is found, we should analyze the possible explanation of the 
trend in terms of time trends in the population measured in different waves 
of such experiments. Basis in this case would first be calculated for the joined 
plane, and then projections of basis vectors of a specific plane provide a 
time dependence of the basis. In all cases we can investigate the difference 
between results obtained for time dependent and approximate time independent 
bases. 

This way a series of mathematical models to deal with the longitudinal 
analysis of biomedical data with repeated categorical measurements and event 
history data can be created. They will be oriented both for population and 
individual prognoses through analysis of individual trajectories in state 
space. For example, if the models are calibrated using NLTCS and Medicare 
data, the results could provide a better understanding of the 
biodemographic mechanisms responsible for aging and morbidity, and could 
be useful in forecasting health and population trends and in stochastic 
investigation in medical economics, including estimating the financial 
effectiveness of new medical technologies.

After such investigation of missing data, they can be filled by probabilities 
$\beta _{jl}^i $, which can be calculated using LLS scores of the known part of 
the outcome pattern and the found linear subspace $\Lambda $. This filling 
probability does not depend on basis choice and on prior information used 
for the basis selection procedure. The largely model 
independent imputation procedure could be very useful in filling categorical data in historical 
cohorts or in dataset where data collection is difficult or costly. 

LLS can be used to analyze data where a high dimensional measurement 
vector represents a hidden structure affecting the results of measurements. 
An abundance of such data recently appeared in genetic studies of 
biological systems where the expression of thousands of genes in cells 
taken from different organs and tissues of humans or laboratory animals is 
measured. Such measurements are performed to find appropriate regularities 
in biological functioning of organs and tissues of respective organisms and to 
detect changes in hidden biological structure due to disease, exposure to 
external factors, aging related changes, etc. Such an analysis will help us 
to better understand mechanisms of genetic regulation, by suggesting genes 
playing key roles in organizing response to changes in internal or external 
milieu.

Forthcoming steps with LLS and the algorithms will include further 
investigation of algorithm properties and development of computer code in 
a form convenient for use by applied researchers. The most important 
property to be investigated is numerical stability of the algorithm. 
Different possibilities to perform subtasks should be investigated to avoid 
ill-posed problems. Several steps in algorithms are based on heuristic 
approaches, so they have to be theoretically investigated and improved. This 
will, on the one hand, provide a basis for selection of numerical methods 
used in the algorithm and, on the other hand, provide rigorous proofs of 
conditions of applicability of the algorithm and reliability of its results. 

\textbf{ACKNOWLEDGEMENTS }: This research was supported by grant P01 
AG17937-05 from National Institute on Aging.

\section*{Appendix}

The mixing distribution (which is a distribution in Euclidean space $R^{|L|}$
carried by a polyhedron defined by (\ref{eq2a})) can be given in form of either
a probabilistic measure $\mu(\beta)$, or probability density function
$\rho_\beta(\beta)$ (or cumulative distribution function $F(\beta)$) with
respect to the standard Lebesque measure in $R^{|L|}$). In the last case one
must allow $\rho_\beta(\beta)$ to be a generalized function to capture
cases where $\mu_\beta$ has a singular component with respect to a Lebesque
measure. 

The permitted domain (\ref{eq2a}) of density function 
$\rho_\beta (\beta )$ in (\ref{eq5}) 
can be reflected explicitly in terms of 
products of $\delta -$ and $\theta -$functions, so the density function has a
representation 
\begin{equation}
\label{eq5a}
\rho_\beta(\beta)=\bar{\rho} (\beta) 
\left(\prod\nolimits_{jl} \theta (\beta _{jl} \right)
\prod\nolimits_j \delta \left( \sum\nolimits_l \beta 
_{jl} -1 \right) 
\end{equation}
Recall that the generalized p.d.f. $\delta(x)$ corresponds to a distribution
concentrated at a single point $x=0$ and is characterized by properties
$\delta(x)=0$ for $x\neq 0$ and $\int\nolimits_{-\infty}^{\infty} dx
\delta(x)=1$;
$\theta(x)$ is a characteristic function 
of positive half-axis, i.e. $\theta(x)=0$
if $x<0$ and  $\theta(x)=1$ if $x\geq 0$. Product 
$\prod_{jl} \theta \prod_j \delta$ is 0 outside the polyhedron $P_\beta$ defined
by conditions (\ref{eq2a}), thus values of $\bar p(\beta)$ outside $P_\beta$ are
insignificant.

In transferring from the integration over $\beta$ to the 
integration over $g$, 
the integral dimension is reduced 
from $\vert L\vert $ with $J$ restrictions defined by $\delta -$functions 
to $K$ with one restriction. Explicitly it can be obtained from (\ref{eq5}) by 
i) orthogonal transformation of integration variables such that $K$ new variables
belong to $\Lambda $; ii) using the fact (main LLS assumption) 
that the density function depends only on these $K$ variables; iii) integration over
remaining $|L|-K$ variables. Formally such a procedure corresponds to
representation of $\bar{\rho}(\beta)$,
\begin{equation}\label{eq7a}
\bar{\rho} (\beta)=\int d{\bf g} \hat{\rho}({\bf g}) \prod\nolimits_{jl} 
\delta \left( \beta _{jl} -\sum\nolimits_k g_k \lambda _{jl}^k \right)
\end{equation} 
and forthcoming integration over all $\beta_{jl}$ using $|L|$ delta functions in
(\ref{eq7a}). $\delta$- and $\theta$-functions extracting permitted integration
area (\ref{eq2a}) in integrals (\ref{eq5}) and (\ref{eq5a}) are transformed such that
permitted area in integration over $g$ are given by (\ref{eq2b}).
   As a result we obtain
unconditional moments $M_\ell$ (and, therefore probabilities $p_\ell$) 
as integrals over $g$ coordinates:
\begin{eqnarray}
\label{eq8}
M_\ell=p_\ell&=&\int d{\bf g}
\left[\prod\nolimits_{j:\ell _j \ne 0} \sum\nolimits_k g_k \lambda^k_{j\ell _j }
\right] 
 \times
\\ && \quad \times
\hat{\rho} ({\bf g})\delta \left( 
\sum\nolimits_k g_k  -1 \right)\,\, \prod\nolimits_{jl} \theta \left( 
\sum\nolimits_k g_k \lambda _{jl}^k  \right) \nonumber
\end{eqnarray}

\textbf{References}

Agresti A (2002). \textit{Categorical Data Analysis}. John Wiley {\&} Sons, Inc. New Jersey.

Akushevich I, Kulminski A, and Manton K (2005) Life tables with covariates: 
Life tables with covariates: Dynamic Model for Nonlinear Analysis of 
Longitudinal Data. \textit{Mathematical Population Studies}, vol. 12(\ref{eq2}), pp.51-80.

Attetkov AV, Galkin SV, and Zarubin VS, Optimization Methods: University 
textbook. N. E. Bauman Moscow State Techn. Univ. Press, Moscow, 2001 (in 
Russian).

Bartholomew, DJ (2002). Old and new approaches to latent variable modeling. 
In: Marcoulides, G. A. and Moustaki, I. (Eds.), \textit{Latent Variable and Latent Structure Models}, (pp.1-14). Mahwah, NJ: 
Lawrence Erlbaum Associates.

Clogg, CC. (1995). \textit{Latent Class Models. }In ``Handbook of Statistical Modeling for the Social and 
Behavioral Sciences'', Arminger, G., Clogg, C.C., Sobel, M.E., eds., New 
York: Plenum Press, 311--360.

Erosheva EA (2002). Grade of Membership and Latent Structure Models with 
Application to Disability Survey Data. Ph.D. Thesis, Department of 
Statistics, Carnegie Mellon University, Pittsburgh, PA, August 2002.

Erosheva EA (2005). Comparing Latent Structures of the Grade of Membership, 
Rasch, and Latent Class Models.~ \textit{Psychometrika (to appear)}.~ 

Forsythe GE, Malcolm MA, and Moler CB. (1977). Computer Methods for 
Mathematical Computations. Prentice-Hall, Inc., Englewood 
Cli{\textregistered}s, NJ.

Goodman, L. A. (1974) Exploratory latent structure analysis using both 
identifiable and unidentifiable models. \textit{Biometrika}, 61, 215-231

Golub GH, Van Loan CF. (1989) Matrix computations. Baltimore :~Johns Hopkins 
University Press.

Haberman SJ (1995). Book review of \textit{Statistical Applications Using Fuzzy Sets}, by Kenneth G. Manton, Max A. Woodbury, 
and H. Dennis Tolley. \textit{Journal of the American Statistical Association}, \textbf{90} (431), 1995, pp. 1131--1133.

Hagenaars J, McCutcheon A (Eds) (2002). \textit{Applied Latent Class Analysis}. Cambridge University Press.

Heinen, T. (1996). \textit{Latent class and discrete latent trait models: Similarities and differences}. Thousand Oaks, California: Sage. 

Kahaner D, Moler C, and Nash S \textit{Numerical Methods and Software} Prentice Hall 1988. 

Kendall MG. and Stuart A (1977). Design and Analysis, and Time-Series, 
volume 3 of The Advanced Theory of Statistics. C. Gri$\pm $n, London, 4 
edition.

Kovtun M, Akushevich I, Manton KG and Tolley HD (2005a). Grade of Membership 
Analysis: One Possible Approach to Foundations. In\textit{ Focus on Probability Theory}, Nova Science Publishers, 
NY. Available at e-print archive arXiv.org, ref. number math.PR/0403373.

Kovtun M, Akushevich I, Manton KG and Tolley HD (2005b). Linear Latent 
Structure Analysis: Mixture Distribution Models with Linear Constrains. 
Submitted to Statistical Methodology. Available at e-print archive 
arXiv.org, ref. number math.PR/.0507025

Kovtun M, Akushevich I, Manton KG and Tolley HD (2005c). A New Efficient 
Algorithm for Construction of LLS Models. Presentation (Poster Session 6: 
Applied Demography, Methods, Health and Mortality) in the Population 
Association of America Annual Meeting, Philadelphia, March 31-April 2. 2005. 
Available at e-print archive arXiv.org, ref. number math.PR/0507021.

Kovtun M, Yashin A, and Akushevich I. (2005d) Convergence of estimators in 
LLS analysis. Submitted for publication to Journal of Nonparametric 
Statistics." Available from e-Print archive arXiv.org at 
http://www.arxiv.org, arXiv code math.PR/0508297, 2005d.

Lazarsfeld PF, Henry NW. (1968) \textit{Latent structure analysis}. Boston: Houghton Mifflin.

Manton KG, Stallard E, Singer BH (1992). Projecting the future size and 
health status of the U.S. elderly population. \textit{International Journal of Forecasting.} 8: P. 433-458.

Manton, KG, Woodbury MA, {\&} Tolley HD (1994). \textit{Statistical applications using fuzzy sets. }John Wiley and Sons, New 
York.

Manton KG, Gu X. (2001) 
Changes in the prevalence of chronic disability in the United States black 
and nonblack population above age 65 from 1982 to 1999. 
\textit{Proceedings of the National Academy of Sciences} 98(11):6354-9.

Manton KG, Lowrimore G, Yashin AI, Kovtun M. (2003) Cluster Analysis of 
Variables. An article for \textit{Encyclopedia of Behavioral Statistics} (B.Evritt, D.Howell, S.Landau, Eds.,).

Mislevy, R. J. (1984). Estimating latent distributions. Psychometrika 49, 
359-381.

Powers DA and Xie Y. (2000) Statistical Methods for Categorical Data 
Analysis. New York: Academic Press. 

Qu Y., Tan M., {\&} Kutner M. H. (1996). Random effects models in latent 
class analysis for evaluating accuracy of diagnostic tests. Biometrics, 52, 
797-810.

Tolley HD, Manton KG. (1992) Large sample properties of estimates of a 
discrete grade of membership model. \textit{Ann Inst Stat Math}. 44(\ref{eq1}):85-95

Uebersax, J. S., {\&} Grove, W. M. (1993). A latent trait finite mixture 
model for the analysis of rating agreement. Biometrics, 49, 823-835.

Uebersax, J. S. (1997). Analysis of student problem behaviors with latent 
trait, latent class, and related probit mixture models. In: Rost J, 
Langeheine R, eds. Applications of Latent Trait and Latent Class Models in 
the Social Sciences. New York, NY: Waxmann; 1997:188-195.

Witteman JC, Grobbee DE, Hofman A (1994) Relation between aortic 
atherosclerosis and blood pressure.Lancet. 25;343, p504-507.

Woodbury MA and Clive J. (1974) Clinical pure types as a fuzzy partition. 
Journal of Cybernetics, 4: 111-121.

Yashin AI, Manton KG. (1997). Effects of unobserved and partially observed 
covariate processes on system failure: A review of models and estimation 
strategies. Statistical Science 12(\ref{eq1}):20-34.

\end{document}